\documentclass{article}
\long\def\ig#1{\relax}
\ig{Thanks to Roberto Minio for this def'n.  Compare the def'n of
\comment in AMSTeX.}

\newcount \coefa
\newcount \coefb
\newcount \coefc
\newcount\tempcounta
\newcount\tempcountb
\newcount\tempcountc
\newcount\tempcountd
\newcount\xext
\newcount\yext
\newcount\xoff
\newcount\yoff
\newcount\gap%
\newcount\arrowtypea
\newcount\arrowtypeb
\newcount\arrowtypec
\newcount\arrowtyped
\newcount\arrowtypee
\newcount\height
\newcount\width
\newcount\xpos
\newcount\ypos
\newcount\run
\newcount\rise
\newcount\arrowlength
\newcount\halflength
\newcount\arrowtype
\newdimen\tempdimen
\newdimen\xlen
\newdimen\ylen
\newsavebox{\tempboxa}%
\newsavebox{\tempboxb}%
\newsavebox{\tempboxc}%

\makeatletter
\def\settypes(#1,#2,#3){\arrowtypea#1 \arrowtypeb#2 \arrowtypec#3}
\def\settoheight#1#2{\setbox\@tempboxa\hbox{#2}#1\ht\@tempboxa\relax}%
\def\settodepth#1#2{\setbox\@tempboxa\hbox{#2}#1\dp\@tempboxa\relax}%
\def\settokens[#1`#2`#3`#4]{%
     \def\tokena{#1}\def\tokenb{#2}\def\tokenc{#3}\def\tokend{#4}}
\def\setsqparms[#1`#2`#3`#4;#5`#6]{%
\arrowtypea #1
\arrowtypeb #2
\arrowtypec #3
\arrowtyped #4
\width #5
\height #6
}
\def\setpos(#1,#2){\xpos=#1 \ypos#2}

\def\bfig{\begin{picture}(\xext,\yext)(\xoff,\yoff)}
\def\efig{\end{picture}}

\def\putbox(#1,#2)#3{\put(#1,#2){\makebox(0,0){$#3$}}}

\def\settriparms[#1`#2`#3;#4]{\settripairparms[#1`#2`#3`1`1;#4]}%

\def\settripairparms[#1`#2`#3`#4`#5;#6]{%
\arrowtypea #1
\arrowtypeb #2
\arrowtypec #3
\arrowtyped #4
\arrowtypee #5
\width #6
\height #6
}

\def\resetparms{\settripairparms[1`1`1`1`1;500]\width 500}%default values%

\resetparms

\def\mvector(#1,#2)#3{%%
\put(0,0){\vector(#1,#2){#3}}%
\put(0,0){\vector(#1,#2){30}}%
}
\def\evector(#1,#2)#3{{%%
\arrowlength #3
\put(0,0){\vector(#1,#2){\arrowlength}}%
\advance \arrowlength by-30
\put(0,0){\vector(#1,#2){\arrowlength}}%
}}

\def\horsize#1#2{%
\settowidth{\tempdimen}{$#2$}%
#1=\tempdimen
\divide #1 by\unitlength
}

\def\vertsize#1#2{%
\settoheight{\tempdimen}{$#2$}%
#1=\tempdimen
\settodepth{\tempdimen}{$#2$}%
\advance #1 by\tempdimen
\divide #1 by\unitlength
}

\def\vertadjust[#1`#2`#3]{%
\vertsize{\tempcounta}{#1}%
\vertsize{\tempcountb}{#2}%
\ifnum \tempcounta<\tempcountb \tempcounta=\tempcountb \fi
\divide\tempcounta by2
\vertsize{\tempcountb}{#3}%
\ifnum \tempcountb>0 \advance \tempcountb by20 \fi
\ifnum \tempcounta<\tempcountb \tempcounta=\tempcountb \fi
}

\def\horadjust[#1`#2`#3]{%
\horsize{\tempcounta}{#1}%
\horsize{\tempcountb}{#2}%
\ifnum \tempcounta<\tempcountb \tempcounta=\tempcountb \fi
\divide\tempcounta by20
\horsize{\tempcountb}{#3}%
\ifnum \tempcountb>0 \advance \tempcountb by60 \fi
\ifnum \tempcounta<\tempcountb \tempcounta=\tempcountb \fi
}

\ig{ In this procedure, #1 is the paramater that sticks out all the way,
#2 sticks out the least and #3 is a label sticking out half way.  #4 is
the amount of the offset.}

\def\sladjust[#1`#2`#3]#4{%
\tempcountc=#4
\horsize{\tempcounta}{#1}%
\divide \tempcounta by2
\horsize{\tempcountb}{#2}%
\divide \tempcountb by2
\advance \tempcountb by-\tempcountc
\ifnum \tempcounta<\tempcountb \tempcounta=\tempcountb\fi
\divide \tempcountc by2
\horsize{\tempcountb}{#3}%
\advance \tempcountb by-\tempcountc
\ifnum \tempcountb>0 \advance \tempcountb by80\fi
\ifnum \tempcounta<\tempcountb \tempcounta=\tempcountb\fi
\advance\tempcounta by20
}

\def\putvector(#1,#2)(#3,#4)#5#6{{%
\xpos=#1
\ypos=#2
\run=#3
\rise=#4
\arrowlength=#5
\arrowtype=#6
\ifnum \arrowtype<0
    \ifnum \run=0
        \advance \ypos by-\arrowlength
    \else
        \tempcounta \arrowlength
        \multiply \tempcounta by\rise
        \divide \tempcounta by\run
        \ifnum\run>0
            \advance \xpos by\arrowlength
            \advance \ypos by\tempcounta
        \else
            \advance \xpos by-\arrowlength
            \advance \ypos by-\tempcounta
        \fi
    \fi
    \multiply \arrowtype by-1
    \multiply \rise by-1
    \multiply \run by-1
\fi
\ifnum \arrowtype=1
    \put(\xpos,\ypos){\vector(\run,\rise){\arrowlength}}%
\else\ifnum \arrowtype=2
    \put(\xpos,\ypos){\mvector(\run,\rise)\arrowlength}%
\else\ifnum\arrowtype=3
    \put(\xpos,\ypos){\evector(\run,\rise){\arrowlength}}%
\fi\fi\fi
}}

\def\putsplitvector(#1,#2)#3#4{%%
\xpos #1
\ypos #2
\arrowtype #4
\halflength #3
\arrowlength #3
\gap 140
\advance \halflength by-\gap
\divide \halflength by2
\ifnum \arrowtype=1
    \put(\xpos,\ypos){\line(0,-1){\halflength}}%
    \advance\ypos by-\halflength
    \advance\ypos by-\gap
    \put(\xpos,\ypos){\vector(0,-1){\halflength}}%
\else\ifnum \arrowtype=2
    \put(\xpos,\ypos){\line(0,-1)\halflength}%
    \put(\xpos,\ypos){\vector(0,-1)3}%
    \advance\ypos by-\halflength
    \advance\ypos by-\gap
    \put(\xpos,\ypos){\vector(0,-1){\halflength}}%
\else\ifnum\arrowtype=3
    \put(\xpos,\ypos){\line(0,-1)\halflength}%
    \advance\ypos by-\halflength
    \advance\ypos by-\gap
    \put(\xpos,\ypos){\evector(0,-1){\halflength}}%
\else\ifnum \arrowtype=-1
    \advance \ypos by-\arrowlength
    \put(\xpos,\ypos){\line(0,1){\halflength}}%
    \advance\ypos by\halflength
    \advance\ypos by\gap
    \put(\xpos,\ypos){\vector(0,1){\halflength}}%
\else\ifnum \arrowtype=-2
    \advance \ypos by-\arrowlength
    \put(\xpos,\ypos){\line(0,1)\halflength}%
    \put(\xpos,\ypos){\vector(0,1)3}%
    \advance\ypos by\halflength
    \advance\ypos by\gap
    \put(\xpos,\ypos){\vector(0,1){\halflength}}%
\else\ifnum\arrowtype=-3
    \advance \ypos by-\arrowlength
    \put(\xpos,\ypos){\line(0,1)\halflength}%
    \advance\ypos by\halflength
    \advance\ypos by\gap
    \put(\xpos,\ypos){\evector(0,1){\halflength}}%
\fi\fi\fi\fi\fi\fi
}

\def\putmorphism(#1)(#2,#3)[#4`#5`#6]#7#8#9{{%
\run #2
\rise #3
\ifnum\rise=0
  \puthmorphism(#1)[#4`#5`#6]{#7}{#8}{#9}%
\else\ifnum\run=0
  \putvmorphism(#1)[#4`#5`#6]{#7}{#8}{#9}%
\else
\setpos(#1)%
\arrowlength #7
\arrowtype #8
\ifnum\run=0
\else\ifnum\rise=0
\else
\ifnum\run>0
    \coefa=1
\else
   \coefa=-1
\fi
\ifnum\arrowtype>0
   \coefb=0
   \coefc=-1
\else
   \coefb=\coefa
   \coefc=1
   \arrowtype=-\arrowtype
\fi
\width=2
\multiply \width by\run
\divide \width by\rise
\ifnum \width<0  \width=-\width\fi
\advance\width by60
\if l#9 \width=-\width\fi
\putbox(\xpos,\ypos){#4}%            %node 1
{\multiply \coefa by\arrowlength%      %node 2
\advance\xpos by\coefa
\multiply \coefa by\rise
\divide \coefa by\run
\advance \ypos by\coefa
\putbox(\xpos,\ypos){#5} }%
{\multiply \coefa by\arrowlength%      %label
\divide \coefa by2
\advance \xpos by\coefa
\advance \xpos by\width
\multiply \coefa by\rise
\divide \coefa by\run
\advance \ypos by\coefa
\if l#9%
   \put(\xpos,\ypos){\makebox(0,0)[r]{$#6$}}%
\else\if r#9%
   \put(\xpos,\ypos){\makebox(0,0)[l]{$#6$}}%
\fi\fi }%
{\multiply \rise by-\coefc%             %arrow
\multiply \run by-\coefc
\multiply \coefb by\arrowlength
\advance \xpos by\coefb
\multiply \coefb by\rise
\divide \coefb by\run
\advance \ypos by\coefb
\multiply \coefc by70
\advance \ypos by\coefc
\multiply \coefc by\run
\divide \coefc by\rise
\advance \xpos by\coefc
\multiply \coefa by140
\multiply \coefa by\run
\divide \coefa by\rise
\advance \arrowlength by\coefa
\ifnum \arrowtype=1
   \put(\xpos,\ypos){\vector(\run,\rise){\arrowlength}}%
\else\ifnum\arrowtype=2
   \put(\xpos,\ypos){\mvector(\run,\rise){\arrowlength}}%
\else\ifnum\arrowtype=3
   \put(\xpos,\ypos){\evector(\run,\rise){\arrowlength}}%
\fi\fi\fi}\fi\fi\fi\fi}}

\def\puthmorphism(#1,#2)[#3`#4`#5]#6#7#8{{%
\xpos #1
\ypos #2
\width #6
\arrowlength #6
\putbox(\xpos,\ypos){#3\vphantom{#4}}%
{\advance \xpos by\arrowlength
\putbox(\xpos,\ypos){\vphantom{#3}#4}}%
\horsize{\tempcounta}{#3}%
\horsize{\tempcountb}{#4}%
\divide \tempcounta by2
\divide \tempcountb by2
\advance \tempcounta by30
\advance \tempcountb by30
\advance \xpos by\tempcounta
\advance \arrowlength by-\tempcounta
\advance \arrowlength by-\tempcountb
\putvector(\xpos,\ypos)(1,0){\arrowlength}{#7}%
\divide \arrowlength by2
\advance \xpos by\arrowlength
\vertsize{\tempcounta}{#5}%
\divide\tempcounta by2
\advance \tempcounta by20
\if a#8 %
   \advance \ypos by\tempcounta
   \putbox(\xpos,\ypos){#5}%
\else
   \advance \ypos by-\tempcounta
   \putbox(\xpos,\ypos){#5}%
\fi}}

\def\putvmorphism(#1,#2)[#3`#4`#5]#6#7#8{{%
\xpos #1
\ypos #2
\arrowlength #6
\arrowtype #7
\settowidth{\xlen}{$#5$}%
\putbox(\xpos,\ypos){#3}%
{\advance \ypos by-\arrowlength
\putbox(\xpos,\ypos){#4}}%
{\advance\arrowlength by-140
\advance \ypos by-70
\ifdim\xlen>0pt
   \if m#8%
      \putsplitvector(\xpos,\ypos){\arrowlength}{\arrowtype}%
   \else
      \putvector(\xpos,\ypos)(0,-1){\arrowlength}{\arrowtype}%
   \fi
\else
   \putvector(\xpos,\ypos)(0,-1){\arrowlength}{\arrowtype}%
\fi}%
\ifdim\xlen>0pt
   \divide \arrowlength by2
   \advance\ypos by-\arrowlength
   \if l#8%
      \advance \xpos by-40
      \put(\xpos,\ypos){\makebox(0,0)[r]{$#5$}}%
   \else\if r#8%
      \advance \xpos by40
      \put(\xpos,\ypos){\makebox(0,0)[l]{$#5$}}%
   \else
      \putbox(\xpos,\ypos){#5}%
   \fi\fi
\fi
}}

\def\topadjust[#1`#2`#3]{%
\yoff=10
\vertadjust[#1`#2`{#3}]%
\advance \yext by\tempcounta
\advance \yext by 10
}
\def\botadjust[#1`#2`#3]{%
\vertadjust[#1`#2`{#3}]%
\advance \yext by\tempcounta
\advance \yoff by-\tempcounta
}
\def\leftadjust[#1`#2`#3]{%
\xoff=0
\horadjust[#1`#2`{#3}]%
\advance \xext by\tempcounta
\advance \xoff by-\tempcounta
}
\def\rightadjust[#1`#2`#3]{%
\horadjust[#1`#2`{#3}]%
\advance \xext by\tempcounta
}
\def\rightsladjust[#1`#2`#3]{%
\sladjust[#1`#2`{#3}]{\width}%
\advance \xext by\tempcounta
}
\def\leftsladjust[#1`#2`#3]{%
\xoff=0
\sladjust[#1`#2`{#3}]{\width}%
\advance \xext by\tempcounta
\advance \xoff by-\tempcounta
}
\def\adjust[#1`#2;#3`#4;#5`#6;#7`#8]{%
\topadjust[#1``{#2}]
\leftadjust[#3``{#4}]
\rightadjust[#5``{#6}]
\botadjust[#7``{#8}]}

\def\putsquarep<#1>(#2)[#3;#4`#5`#6`#7]{{%
\setsqparms[#1]%
\setpos(#2)%
\settokens[#3]%
\puthmorphism(\xpos,\ypos)[\tokenc`\tokend`{#7}]{\width}{\arrowtyped}b%
\advance\ypos by \height
\puthmorphism(\xpos,\ypos)[\tokena`\tokenb`{#4}]{\width}{\arrowtypea}a%
\putvmorphism(\xpos,\ypos)[``{#5}]{\height}{\arrowtypeb}l%
\advance\xpos by \width
\putvmorphism(\xpos,\ypos)[``{#6}]{\height}{\arrowtypec}r%
}}

\def\putsquare{\@ifnextchar <{\putsquarep}{\putsquarep%
   <\arrowtypea`\arrowtypeb`\arrowtypec`\arrowtyped;\width`\height>}}
\def\square{\@ifnextchar< {\squarep}{\squarep
   <\arrowtypea`\arrowtypeb`\arrowtypec`\arrowtyped;\width`\height>}}
\def\squarep<#1>[#2`#3`#4`#5;#6`#7`#8`#9]{{%          %     #2------>#3
\setsqparms[#1]%                                      %      |       |
\xext=\width                                          %      |       |
\yext=\height                                         %    #7|       |#8
\topadjust[#2`#3`{#6}]%                               %      |       |
\botadjust[#4`#5`{#9}]%                               %      |       |
\leftadjust[#2`#4`{#7}]%                              %
\rightadjust[#3`#5`{#8}]%                             %     #4------>#5
\begin{picture}(\xext,\yext)(\xoff,\yoff)%                      #9
\putsquarep<\arrowtypea`\arrowtypeb`\arrowtypec`\arrowtyped;\width`\height>%
(0,0)[#2`#3`#4`#5;#6`#7`#8`{#9}]%
\end{picture}%
}}

\def\putptrianglep<#1>(#2,#3)[#4`#5`#6;#7`#8`#9]{{%
\settriparms[#1]%
\xpos=#2 \ypos=#3
\advance\ypos by \height
\puthmorphism(\xpos,\ypos)[#4`#5`{#7}]{\height}{\arrowtypea}a%
\putvmorphism(\xpos,\ypos)[`#6`{#8}]{\height}{\arrowtypeb}l%
\advance\xpos by\height
\putmorphism(\xpos,\ypos)(-1,-1)[``{#9}]{\height}{\arrowtypec}r%
}}

\def\putptriangle{\@ifnextchar <{\putptrianglep}{\putptrianglep
   <\arrowtypea`\arrowtypeb`\arrowtypec;\height>}}
\def\ptriangle{\@ifnextchar <{\ptrianglep}{\ptrianglep
   <\arrowtypea`\arrowtypeb`\arrowtypec;\height>}}

\def\ptrianglep<#1>[#2`#3`#4;#5`#6`#7]{{%%       #5
\settriparms[#1]%
\width=\height                         %      #2----->#3
\xext=\width                           %      |      /
\yext=\width                           %      |     /
\topadjust[#2`#3`{#5}]%                %    #6|    /#7
\botadjust[#3``]%                      %      |   /
\leftadjust[#2`#4`{#6}]%               %      |  /
\rightsladjust[#3`#4`{#7}]%            %
\begin{picture}(\xext,\yext)(\xoff,\yoff)%    #4
\putptrianglep<\arrowtypea`\arrowtypeb`\arrowtypec;\height>%
(0,0)[#2`#3`#4;#5`#6`{#7}]%
\end{picture}%
}}

\def\putqtrianglep<#1>(#2,#3)[#4`#5`#6;#7`#8`#9]{{%
\settriparms[#1]%
\xpos=#2 \ypos=#3
\advance\ypos by\height
\puthmorphism(\xpos,\ypos)[#4`#5`{#7}]{\height}{\arrowtypea}a%
\putmorphism(\xpos,\ypos)(1,-1)[``{#8}]{\height}{\arrowtypeb}l%
\advance\xpos by\height
\putvmorphism(\xpos,\ypos)[`#6`{#9}]{\height}{\arrowtypec}r%
}}

\def\putqtriangle{\@ifnextchar <{\putqtrianglep}{\putqtrianglep
   <\arrowtypea`\arrowtypeb`\arrowtypec;\height>}}
\def\qtriangle{\@ifnextchar <{\qtrianglep}{\qtrianglep
   <\arrowtypea`\arrowtypeb`\arrowtypec;\height>}}

\def\qtrianglep<#1>[#2`#3`#4;#5`#6`#7]{{%%
\settriparms[#1]%                                  #5
\width=\height                         %        #2----->#3
\xext=\width                           %         \      |
\yext=\height                          %          \     |
\topadjust[#2`#3`{#5}]%                %         #6\    |#7
\botadjust[#4``]%                      %            \   |
\leftsladjust[#2`#4`{#6}]%             %             \  |
\rightadjust[#3`#4`{#7}]%              %
\begin{picture}(\xext,\yext)(\xoff,\yoff)%             #4
\putqtrianglep<\arrowtypea`\arrowtypeb`\arrowtypec;\height>%
(0,0)[#2`#3`#4;#5`#6`{#7}]%
\end{picture}%
}}

\def\putdtrianglep<#1>(#2,#3)[#4`#5`#6;#7`#8`#9]{{%
\settriparms[#1]%
\xpos=#2 \ypos=#3
\puthmorphism(\xpos,\ypos)[#5`#6`{#9}]{\height}{\arrowtypec}b%
\advance\xpos by \height \advance\ypos by\height
\putmorphism(\xpos,\ypos)(-1,-1)[``{#7}]{\height}{\arrowtypea}l%
\putvmorphism(\xpos,\ypos)[#4``{#8}]{\height}{\arrowtypeb}r%
}}

\def\putdtriangle{\@ifnextchar <{\putdtrianglep}{\putdtrianglep
   <\arrowtypea`\arrowtypeb`\arrowtypec;\height>}}
\def\dtriangle{\@ifnextchar <{\dtrianglep}{\dtrianglep
   <\arrowtypea`\arrowtypeb`\arrowtypec;\height>}}

\def\dtrianglep<#1>[#2`#3`#4;#5`#6`#7]{{%%
\settriparms[#1]%                                          #2
\width=\height                         %                  / |
\xext=\width                           %                 /  |
\yext=\height                          %              #5/   |#6
\topadjust[#2``]%                      %               /    |
\botadjust[#3`#4`{#7}]%                %              /     |
\leftsladjust[#3`#2`{#5}]%             %
\rightadjust[#2`#4`{#6}]%              %            #3----->#4
\begin{picture}(\xext,\yext)(\xoff,\yoff)%              #7
\putdtrianglep<\arrowtypea`\arrowtypeb`\arrowtypec;\height>%
(0,0)[#2`#3`#4;#5`#6`{#7}]%
\end{picture}%
}}

\def\putbtrianglep<#1>(#2,#3)[#4`#5`#6;#7`#8`#9]{{%
\settriparms[#1]%
\xpos=#2 \ypos=#3
\puthmorphism(\xpos,\ypos)[#5`#6`{#9}]{\height}{\arrowtypec}b%
\advance\ypos by\height
\putmorphism(\xpos,\ypos)(1,-1)[``{#8}]{\height}{\arrowtypeb}r%
\putvmorphism(\xpos,\ypos)[#4``{#7}]{\height}{\arrowtypea}l%
}}

\def\putbtriangle{\@ifnextchar <{\putbtrianglep}{\putbtrianglep
   <\arrowtypea`\arrowtypeb`\arrowtypec;\height>}}
\def\btriangle{\@ifnextchar <{\btrianglep}{\btrianglep
   <\arrowtypea`\arrowtypeb`\arrowtypec;\height>}}

\def\btrianglep<#1>[#2`#3`#4;#5`#6`#7]{{%%
\settriparms[#1]%                                     #2
\width=\height                         %              | \
\xext=\width                           %              |  \
\yext=\height                          %            #5|   \#6
\topadjust[#2``]%                      %              |    \
\botadjust[#3`#4`{#7}]%                %              |     \
\leftadjust[#2`#3`{#5}]%               %
\rightsladjust[#4`#2`{#6}]%            %              #3----->#4
\begin{picture}(\xext,\yext)(\xoff,\yoff)%                #7
\putbtrianglep<\arrowtypea`\arrowtypeb`\arrowtypec;\height>%
(0,0)[#2`#3`#4;#5`#6`{#7}]%
\end{picture}%
}}

\def\putAtrianglep<#1>(#2,#3)[#4`#5`#6;#7`#8`#9]{{%
\settriparms[#1]%
\xpos=#2 \ypos=#3
{\multiply \height by2
\puthmorphism(\xpos,\ypos)[#5`#6`{#9}]{\height}{\arrowtypec}b}%
\advance\xpos by\height \advance\ypos by\height
\putmorphism(\xpos,\ypos)(-1,-1)[#4``{#7}]{\height}{\arrowtypea}l%
\putmorphism(\xpos,\ypos)(1,-1)[``{#8}]{\height}{\arrowtypeb}r%
}}

\def\putAtriangle{\@ifnextchar <{\putAtrianglep}{\putAtrianglep
   <\arrowtypea`\arrowtypeb`\arrowtypec;\height>}}
\def\Atriangle{\@ifnextchar <{\Atrianglep}{\Atrianglep
   <\arrowtypea`\arrowtypeb`\arrowtypec;\height>}}

\def\Atrianglep<#1>[#2`#3`#4;#5`#6`#7]{{%%
\settriparms[#1]%                                 #2
\width=\height                         %         /   \
\xext=\width                           %        /     \
\yext=\height                          %     #5/       \#6
\topadjust[#2``]%                      %      /         \
\botadjust[#3`#4`{#7}]%                %     /           \
\multiply \xext by2 %                  %
\leftsladjust[#3`#2`{#5}]%             %   #3------------>#4
\rightsladjust[#4`#2`{#6}]%            %          #7
\begin{picture}(\xext,\yext)(\xoff,\yoff)%
\putAtrianglep<\arrowtypea`\arrowtypeb`\arrowtypec;\height>%
(0,0)[#2`#3`#4;#5`#6`{#7}]%
\end{picture}%
}}

\def\putAtrianglepairp<#1>(#2)[#3;#4`#5`#6`#7`#8]{{
\settripairparms[#1]%
\setpos(#2)%
\settokens[#3]%
\puthmorphism(\xpos,\ypos)[\tokenb`\tokenc`{#7}]{\height}{\arrowtyped}b%
\advance\xpos by\height
\advance\ypos by\height
\putmorphism(\xpos,\ypos)(-1,-1)[\tokena``{#4}]{\height}{\arrowtypea}l%
\putvmorphism(\xpos,\ypos)[``{#5}]{\height}{\arrowtypeb}m%
\putmorphism(\xpos,\ypos)(1,-1)[``{#6}]{\height}{\arrowtypec}r%
}}

\def\putAtrianglepair{\@ifnextchar <{\putAtrianglepairp}{\putAtrianglepairp%
   <\arrowtypea`\arrowtypeb`\arrowtypec`\arrowtyped`\arrowtypee;\height>}}
\def\Atrianglepair{\@ifnextchar <{\Atrianglepairp}{\Atrianglepairp%
   <\arrowtypea`\arrowtypeb`\arrowtypec`\arrowtyped`\arrowtypee;\height>}}

\def\Atrianglepairp<#1>[#2;#3`#4`#5`#6`#7]{{%
\settripairparms[#1]%
\settokens[#2]%
\width=\height
\xext=\width
\yext=\height
\topadjust[\tokena``]%
\vertadjust[\tokenb`\tokenc`{#6}]%                      %  #2a
\tempcountd=\tempcounta                       %           / | \
\vertadjust[\tokenc`\tokend`{#7}]%            %          /  |  \
\ifnum\tempcounta<\tempcountd                 %       #3/  #4   \#5
\tempcounta=\tempcountd\fi                    %        /    |    \
\advance \yext by\tempcounta                  %       /     |     \
\advance \yoff by-\tempcounta                 %
\multiply \xext by2 %                         %     #2b---->#2c---->#2d
\leftsladjust[\tokenb`\tokena`{#3}]%          %         #6     #7
\rightsladjust[\tokend`\tokena`{#5}]%
\begin{picture}(\xext,\yext)(\xoff,\yoff)%
\putAtrianglepairp
<\arrowtypea`\arrowtypeb`\arrowtypec`\arrowtyped`\arrowtypee;\height>%
(0,0)[#2;#3`#4`#5`#6`{#7}]%
\end{picture}%
}}

\def\putVtrianglep<#1>(#2,#3)[#4`#5`#6;#7`#8`#9]{{%
\settriparms[#1]%
\xpos=#2 \ypos=#3
\advance\ypos by\height
{\multiply\height by2
\puthmorphism(\xpos,\ypos)[#4`#5`{#7}]{\height}{\arrowtypea}a}%
\putmorphism(\xpos,\ypos)(1,-1)[`#6`{#8}]{\height}{\arrowtypeb}l%
\advance\xpos by\height
\advance\xpos by\height
\putmorphism(\xpos,\ypos)(-1,-1)[``{#9}]{\height}{\arrowtypec}r%
}}

\def\putVtriangle{\@ifnextchar <{\putVtrianglep}{\putVtrianglep
   <\arrowtypea`\arrowtypeb`\arrowtypec;\height>}}
\def\Vtriangle{\@ifnextchar <{\Vtrianglep}{\Vtrianglep
   <\arrowtypea`\arrowtypeb`\arrowtypec;\height>}}

\def\Vtrianglep<#1>[#2`#3`#4;#5`#6`#7]{{%%
\settriparms[#1]%                                      #5
\width=\height                         %        #2------------->#3
\xext=\width                           %         \             /
\yext=\height                          %          \           /
\topadjust[#2`#3`{#5}]%                %         #6\         /#7
\botadjust[#4``]%                      %            \       /
\multiply \xext by2 %                  %             \     /
\leftsladjust[#2`#3`{#6}]%             %
\rightsladjust[#3`#4`{#7}]%            %               #4
\begin{picture}(\xext,\yext)(\xoff,\yoff)%
\putVtrianglep<\arrowtypea`\arrowtypeb`\arrowtypec;\height>%
(0,0)[#2`#3`#4;#5`#6`{#7}]%
\end{picture}%
}}

\def\putVtrianglepairp<#1>(#2)[#3;#4`#5`#6`#7`#8]{{
\settripairparms[#1]%
\setpos(#2)%
\settokens[#3]%
\advance\ypos by\height
\putmorphism(\xpos,\ypos)(1,-1)[`\tokend`{#6}]{\height}{\arrowtypec}l%
\puthmorphism(\xpos,\ypos)[\tokena`\tokenb`{#4}]{\height}{\arrowtypea}a%
\advance\xpos by\height
\putvmorphism(\xpos,\ypos)[``{#7}]{\height}{\arrowtyped}m%
\advance\xpos by\height
\putmorphism(\xpos,\ypos)(-1,-1)[``{#8}]{\height}{\arrowtypee}r%
}}

\def\putVtrianglepair{\@ifnextchar <{\putVtrianglepairp}{\putVtrianglepairp%
    <\arrowtypea`\arrowtypeb`\arrowtypec`\arrowtyped`\arrowtypee;\height>}}
\def\Vtrianglepair{\@ifnextchar <{\Vtrianglepairp}{\Vtrianglepairp%
    <\arrowtypea`\arrowtypeb`\arrowtypec`\arrowtyped`\arrowtypee;\height>}}

\def\Vtrianglepairp<#1>[#2;#3`#4`#5`#6`#7]{{%
\settripairparms[#1]%
\settokens[#2]%                            #3      #4
\xext=\height                  %        #2a---->#2b---->#2c
\width=\height                 %         \      |      /
\yext=\height                  %          \     |     /
\vertadjust[\tokena`\tokenb`{#4}]%       #5\   #6    /#7
\tempcountd=\tempcounta        %            \   |   /
\vertadjust[\tokenb`\tokenc`{#5}]%           \  |  /
\ifnum\tempcounta<\tempcountd%
\tempcounta=\tempcountd\fi%                    #2d
\advance \yext by\tempcounta
\botadjust[\tokend``]%
\multiply \xext by2
\leftsladjust[\tokena`\tokend`{#6}]%
\rightsladjust[\tokenc`\tokend`{#7}]%
\begin{picture}(\xext,\yext)(\xoff,\yoff)%
\putVtrianglepairp
<\arrowtypea`\arrowtypeb`\arrowtypec`\arrowtyped`\arrowtypee;\height>%
(0,0)[#2;#3`#4`#5`#6`{#7}]%
\end{picture}%
}}

\def\putCtrianglep<#1>(#2,#3)[#4`#5`#6;#7`#8`#9]{{%
\settriparms[#1]%
\xpos=#2 \ypos=#3
\advance\ypos by\height
\putmorphism(\xpos,\ypos)(1,-1)[``{#9}]{\height}{\arrowtypec}l%
\advance\xpos by\height
\advance\ypos by\height
\putmorphism(\xpos,\ypos)(-1,-1)[#4`#5`{#7}]{\height}{\arrowtypea}l%
{\multiply\height by 2
\putvmorphism(\xpos,\ypos)[`#6`{#8}]{\height}{\arrowtypeb}r}%
}}

\def\putCtriangle{\@ifnextchar <{\putCtrianglep}{\putCtrianglep
    <\arrowtypea`\arrowtypeb`\arrowtypec;\height>}}
\def\Ctriangle{\@ifnextchar <{\Ctrianglep}{\Ctrianglep
    <\arrowtypea`\arrowtypeb`\arrowtypec;\height>}}

\def\Ctrianglep<#1>[#2`#3`#4;#5`#6`#7]{{%%
\settriparms[#1]%                                         #2
\width=\height                          %                / |
\xext=\width                            %               /  |
\yext=\height                           %            #5/   |
\multiply \yext by2 %                   %             /    |
\topadjust[#2``]%                       %            /     |
\botadjust[#4``]%                       %           v      |
\sladjust[#3`#2`{#5}]{\width}%          %          #3      |#6
\tempcountd=\tempcounta                 %           \      |
\sladjust[#3`#4`{#7}]{\width}%          %            \     |
\ifnum \tempcounta<\tempcountd          %           #7\    |
\tempcounta=\tempcountd\fi              %              \   |
\advance \xext by\tempcounta            %               \  |
\advance \xoff by-\tempcounta           %
\rightadjust[#2`#4`{#6}]%               %                 #4
\begin{picture}(\xext,\yext)(\xoff,\yoff)%
\putCtrianglep<\arrowtypea`\arrowtypeb`\arrowtypec;\height>%
(0,0)[#2`#3`#4;#5`#6`{#7}]%
\end{picture}%
}}

\def\putDtrianglep<#1>(#2,#3)[#4`#5`#6;#7`#8`#9]{{%
\settriparms[#1]%
\xpos=#2 \ypos=#3
\advance\xpos by\height \advance\ypos by\height
\putmorphism(\xpos,\ypos)(-1,-1)[``{#9}]{\height}{\arrowtypec}r%
\advance\xpos by-\height \advance\ypos by\height
\putmorphism(\xpos,\ypos)(1,-1)[`#5`{#8}]{\height}{\arrowtypeb}r%
{\multiply\height by 2
\putvmorphism(\xpos,\ypos)[#4`#6`{#7}]{\height}{\arrowtypea}l}%
}}

\def\putDtriangle{\@ifnextchar <{\putDtrianglep}{\putDtrianglep
    <\arrowtypea`\arrowtypeb`\arrowtypec;\height>}}
\def\Dtriangle{\@ifnextchar <{\Dtrianglep}{\Dtrianglep
   <\arrowtypea`\arrowtypeb`\arrowtypec;\height>}}

\def\Dtrianglep<#1>[#2`#3`#4;#5`#6`#7]{{%%
\settriparms[#1]%                                 #2
\width=\height                         %          | \
\xext=\height                          %          |  \
\yext=\height                          %          |   \#6
\multiply \yext by2 %                  %          |    \
\topadjust[#2``]%                      %          |     \
\botadjust[#4``]%                      %          |
\leftadjust[#2`#4`{#5}]%               %        #5|      #3
\sladjust[#3`#2`{#5}]{\height}%        %          |      /
\tempcountd=\tempcountd                %          |     /
\sladjust[#3`#4`{#7}]{\height}%        %          |    /#7
\ifnum \tempcounta<\tempcountd         %          |   /
\tempcounta=\tempcountd\fi             %          |  /
\advance \xext by\tempcounta           %
\begin{picture}(\xext,\yext)(\xoff,\yoff)%        #4
\putDtrianglep<\arrowtypea`\arrowtypeb`\arrowtypec;\height>%
(0,0)[#2`#3`#4;#5`#6`{#7}]%
\end{picture}%
}}

\def\setrecparms[#1`#2]{\width=#1 \height=#2}%
\def\recursep<#1`#2>[#3;#4`#5`#6`#7`#8]{{%
\width=#1 \height=#2
\settokens[#3]
\settowidth{\tempdimen}{$\tokena$}
\ifdim\tempdimen=0pt
  \savebox{\tempboxa}{\hbox{$\tokenb$}}%
  \savebox{\tempboxb}{\hbox{$\tokend$}}%
  \savebox{\tempboxc}{\hbox{$#6$}}%
\else
  \savebox{\tempboxa}{\hbox{$\hbox{$\tokena$}\times\hbox{$\tokenb$}$}}%
  \savebox{\tempboxb}{\hbox{$\hbox{$\tokena$}\times\hbox{$\tokend$}$}}%
  \savebox{\tempboxc}{\hbox{$\hbox{$\tokena$}\times\hbox{$#6$}$}}%
\fi
\ypos=\height
\divide\ypos by 2
\xpos=\ypos
\advance\xpos by \width
\xext=\xpos \yext=\height
\topadjust[#3`\usebox{\tempboxa}`{#4}]%
\botadjust[#5`\usebox{\tempboxb}`{#8}]%
\sladjust[\tokenc`\tokenb`{#5}]{\ypos}%
\tempcountd=\tempcounta
\sladjust[\tokenc`\tokend`{#5}]{\ypos}%
\ifnum \tempcounta<\tempcountd
\tempcounta=\tempcountd\fi
\advance \xext by\tempcounta
\advance \xoff by-\tempcounta
\rightadjust[\usebox{\tempboxa}`\usebox{\tempboxb}`\usebox{\tempboxc}]%
\bfig
\putCtrianglep<-1`1`1;\ypos>(0,0)[`\tokenc`;#5`#6`{#7}]%
\puthmorphism(\ypos,0)[\tokend`\usebox{\tempboxb}`{#8}]{\width}{-1}b%
\puthmorphism(\ypos,\height)[\tokenb`\usebox{\tempboxa}`{#4}]{\width}{-1}a%
\advance\ypos by \width
\putvmorphism(\ypos,\height)[``\usebox{\tempboxc}]{\height}1r%
\efig
}}

\def\recurse{\@ifnextchar <{\recursep}{\recursep<\width`\height>}}

\def\puttwohmorphisms(#1,#2)[#3`#4;#5`#6]#7#8#9{{%
\puthmorphism(#1,#2)[#3`#4`]{#7}0a
\ypos=#2
\advance\ypos by 20
\puthmorphism(#1,\ypos)[\phantom{#3}`\phantom{#4}`#5]{#7}{#8}a
\advance\ypos by -40
\puthmorphism(#1,\ypos)[\phantom{#3}`\phantom{#4}`#6]{#7}{#9}b
}}

\def\puttwovmorphisms(#1,#2)[#3`#4;#5`#6]#7#8#9{{%
\putvmorphism(#1,#2)[#3`#4`]{#7}0a
\xpos=#1
\advance\xpos by -20
\putvmorphism(\xpos,#2)[\phantom{#3}`\phantom{#4}`#5]{#7}{#8}l
\advance\xpos by 40
\putvmorphism(\xpos,#2)[\phantom{#3}`\phantom{#4}`#6]{#7}{#9}r
}}

\def\puthcoequalizer(#1)[#2`#3`#4;#5`#6`#7]#8#9{{%
\setpos(#1)%
\puttwohmorphisms(\xpos,\ypos)[#2`#3;#5`#6]{#8}11%
\advance\xpos by #8
\puthmorphism(\xpos,\ypos)[\phantom{#3}`#4`#7]{#8}1{#9}
}}

\def\putvcoequalizer(#1)[#2`#3`#4;#5`#6`#7]#8#9{{%
\setpos(#1)%
\puttwovmorphisms(\xpos,\ypos)[#2`#3;#5`#6]{#8}11%
\advance\ypos by -#8
\putvmorphism(\xpos,\ypos)[\phantom{#3}`#4`#7]{#8}1{#9}
}}

\def\putthreehmorphisms(#1)[#2`#3;#4`#5`#6]#7(#8)#9{{%
\setpos(#1) \settypes(#8)
\if a#9 %
     \vertsize{\tempcounta}{#5}%
     \vertsize{\tempcountb}{#6}%
     \ifnum \tempcounta<\tempcountb \tempcounta=\tempcountb \fi
\else
     \vertsize{\tempcounta}{#4}%
     \vertsize{\tempcountb}{#5}%
     \ifnum \tempcounta<\tempcountb \tempcounta=\tempcountb \fi
\fi
\advance \tempcounta by 60
\puthmorphism(\xpos,\ypos)[#2`#3`#5]{#7}{\arrowtypeb}{#9}
\advance\ypos by \tempcounta
\puthmorphism(\xpos,\ypos)[\phantom{#2}`\phantom{#3}`#4]{#7}{\arrowtypea}{#9}
\advance\ypos by -\tempcounta \advance\ypos by -\tempcounta
\puthmorphism(\xpos,\ypos)[\phantom{#2}`\phantom{#3}`#6]{#7}{\arrowtypec}{#9}
}}

\def\putarc(#1,#2)[#3`#4`#5]#6#7#8{{%
\xpos #1
\ypos #2
\width #6
\arrowlength #6
\putbox(\xpos,\ypos){#3\vphantom{#4}}%
{\advance \xpos by\arrowlength
\putbox(\xpos,\ypos){\vphantom{#3}#4}}%
\horsize{\tempcounta}{#3}%
\horsize{\tempcountb}{#4}%
\divide \tempcounta by2
\divide \tempcountb by2
\advance \tempcounta by30
\advance \tempcountb by30
\advance \xpos by\tempcounta
\advance \arrowlength by-\tempcounta
\advance \arrowlength by-\tempcountb
\halflength=\arrowlength \divide\halflength by 2
\divide\arrowlength by 5
\put(\xpos,\ypos){\bezier{\arrowlength}(0,0)(50,50)(\halflength,50)}
\ifnum #7=-1 \put(\xpos,\ypos){\vector(-3,-2)0} \fi
\advance\xpos by \halflength
\put(\xpos,\ypos){\xpos=\halflength \advance\xpos by -50
   \bezier{\arrowlength}(0,50)(\xpos,50)(\halflength,0)}
\ifnum #7=1 {\advance \xpos by
   \halflength \put(\xpos,\ypos){\vector(3,-2)0}} \fi
\advance\ypos by 50
\vertsize{\tempcounta}{#5}%
\divide\tempcounta by2
\advance \tempcounta by20
\if a#8 %
   \advance \ypos by\tempcounta
   \putbox(\xpos,\ypos){#5}%
\else
   \advance \ypos by-\tempcounta
   \putbox(\xpos,\ypos){#5}%
\fi
}}

\makeatother

\usepackage{amsthm}
\usepackage{dsfont}
\usepackage{stmaryrd}
\newtheorem{theorem}{Theorem}[section]

\newtheorem{corollary}[theorem]{Corollary}
\newtheorem{fact}[theorem]{Fact}

\newtheorem{proposition}[theorem]{Proposition}

\begin{document}

\sloppy

%commands
\newcommand{\nl}{\hspace{2cm}\\ }

\def\nec{\Box}
\def\pos{\Diamond}
\def\diam{{\tiny\Diamond}}

\def\lc{\lceil}
\def\rc{\rceil}
\def\lf{\lfloor}
\def\rf{\rfloor}
\def\lk{\langle}
\def\rk{\rangle}
\def\blk{\dot{\langle\!\!\langle}}
\def\brk{\dot{\rangle\!\!\rangle}}

\newcommand{\pa}{\parallel}
\newcommand{\lra}{\longrightarrow}
\newcommand{\hra}{\hookrightarrow}
\newcommand{\hla}{\hookleftarrow}
\newcommand{\ra}{\rightarrow}
\newcommand{\la}{\leftarrow}
\newcommand{\lla}{\longleftarrow}
\newcommand{\da}{\downarrow}
\newcommand{\ua}{\uparrow}
\newcommand{\dA}{\downarrow\!\!\!^\bullet}
\newcommand{\uA}{\uparrow\!\!\!_\bullet}
\newcommand{\Da}{\Downarrow}
\newcommand{\DA}{\Downarrow\!\!\!^\bullet}
\newcommand{\UA}{\Uparrow\!\!\!_\bullet}
\newcommand{\Ua}{\Uparrow}
\newcommand{\Lra}{\Longrightarrow}
\newcommand{\Ra}{\Rightarrow}
\newcommand{\Lla}{\Longleftarrow}
\newcommand{\La}{\Leftarrow}
\newcommand{\nperp}{\perp\!\!\!\!\!\setminus\;\;}
\newcommand{\pq}{\preceq}

\newcommand{\lms}{\longmapsto}
\newcommand{\ms}{\mapsto}
\newcommand{\subseteqnot}{\subseteq\hskip-4 mm_\not\hskip3 mm}

\def\o{{\omega}}

\def\bA{{\bf A}}
\def\bEM{{\bf EM}}
\def\bM{{\bf M}}
\def\bN{{\bf N}}
\def\bF{{\bf F}}
\def\bC{{\bf C}}
\def\bI{{\bf I}}
\def\bK{{\bf K}}
\def\bL{{\bf L}}
\def\bT{{\bf T}}
\def\bS{{\bf S}}
\def\bD{{\bf D}}
\def\bB{{\bf B}}
\def\bW{{\bf W}}
\def\bP{{\bf P}}
\def\bX{{\bf X}}
\def\bY{{\bf Y}}
\def\bZ{{\bf Z}}
\def\ba{{\bf a}}
\def\bb{{\bf b}}
\def\bc{{\bf c}}
\def\bd{{\bf d}}
\def\bh{{\bf h}}
\def\bi{{\bf i}}
\def\bj{{\bf j}}
\def\bk{{\bf k}}
\def\bm{{\bf m}}
\def\bn{{\bf n}}
\def\bp{{\bf p}}
\def\bq{{\bf q}}
\def\be{{\bf e}}
\def\br{{\bf r}}
\def\bi{{\bf i}}
\def\bs{{\bf s}}
\def\bt{{\bf t}}
\def\jeden{{\bf 1}}
\def\dwa{{\bf 2}}
\def\trzy{{\bf 3}}

\def\cB{{\cal B}}
\def\cA{{\cal A}}
\def\cC{{\cal C}}
\def\cD{{\cal D}}
\def\cE{{\cal E}}
\def\cEM{{\cal EM}}
\def\cF{{\cal F}}
\def\cG{{\cal G}}
\def\cI{{\cal I}}
\def\cJ{{\cal J}}
\def\cK{{\cal K}}
\def\cL{{\cal L}}
\def\cN{{\cal N}}
\def\cM{{\cal M}}
\def\cO{{\cal O}}
\def\cP{{\cal P}}
\def\cQ{{\cal Q}}
\def\cR{{\cal R}}
\def\cS{{\cal S}}
\def\cT{{\cal T}}
\def\cU{{\cal U}}
\def\cV{{\cal V}}
\def\cW{{\cal W}}
\def\cX{{\cal X}}
\def\cY{{\cal Y}}

%categories

%of functors and monads
\def\Mnd{{\bf Mnd}}
\def\MND{{\bf MND}}
\def\AMnd{{\bf AnMnd}}
\def\ANMND{{\bf ANMND}}
\def\An{{\bf An}}
\def\AN{{\bf AN}}
\def\Poly{{\bf Poly}}
\def\SAN{{\bf SAN}}
\def\San{{\bf San}}
\def\iSanMnd{_\infty{\bf SanMnd}}
\def\iSan{_\infty{\bf San}}
\def\Taut{{\bf Taut}}
\def\PMnd{{\bf PolyMnd}}
\def\SanMnd{{\bf SanMnd}}
\def\SANMND{{\bf SANMND}}
\def\SanMnd{{\bf SanMnd}}
\def\RiMnd{{\bf RiMnd}}
\def\End{{\bf End}}
\def\END{{\bf END}}

%of theories
\def\ET{\bf ET}
\def\RegET{\bf RegET}
\def\RET{\bf RegET}
\def\LrET{\bf LrET}
\def\RiET{\bf RiET}
\def\SregET{\bf SregET}
\def\Cart{\bf Cart}
\def\wCart{\bf wCart}
\def\CartMnd{\bf CartMnd}
\def\wCartMnd{\bf wCartMnd}

%of Lawvere theories
\def\LT{\bf LT}
\def\RegLT{\bf RegLT}
\def\ALT{\bf AnLT}
\def\RiLT{\bf RiLT}

%of Operads
\def\FOp{\bf FOp}
\def\RegOp{\bf RegOp}
\def\SOp{\bf SOp}
\def\RiOp{\bf RiOp}

%various other categories and such
\def\Ord{{\bf Ord}}
\def\Card{{\bf Card}}
\def\CAT{{{\bf CAT}}}
\def\MonCat{{{\bf MonCat}}}
\def\Mon{{{\bf Mon}}}
\def\Cat{{{\bf Cat}}}

\def\F{\mathds{F}}
\def\S{\mathds{S}}
\def\I{\mathds{I}}
\def\B{\mathds{B}}

%functors
\def\V{\mathds{V}}
\def\W{\mathds{W}}
\def\M{\mathds{M}}
\def\N{\mathds{N}}
\def\R{\mathds{R}}

\def\Op{{\cal O}p}

\def\Vb{\bar{\mathds{V}}}
\def\Wb{\bar{\mathds{W}}}

\def\Sym{{\cal S}ym}

%leftovers
\def\P{{\cal P}}
\def\Q{{\cal Q}}

\pagenumbering{arabic} \setcounter{page}{1}

\title{\bf\Large Generalized P\l onka Sums and Products}

\author{ Marek Zawadowski\\
Instytut Matematyki, Uniwersytet Warszawski\\
ul. S.Banacha 2,\\
00-913 Warszawa, Poland\\
zawado@mimuw.edu.pl\\
%\date{October 25, 2012}
}

\maketitle
\begin{flushright}
\em Dedicated to George Janelidze\\ on the occasion of his
60th birthday.
\end{flushright}

\begin{abstract} We give an abstract categorical treatment of P\l onka sums and products using lax and oplax morphisms of monads. P\l onka sums were originally defined as operations on algebras of regular theories. Their arities are sup-semilattices. It turns out that even more general operations are available on the categories of algebras of semi-analytic monads. Their arities are the categories of the regular polynomials over any sup-semilattice, i.e. any algebra for the terminal semi-analytic monad.  We also show that similar operations can be defined on any category of algebras of any analytic monad. This time we can allow the arities to be the categories of linear polynomials over any commutative monoid, i.e. any algebra for the terminal analytic monad.  There are also dual operations of P\l onka products. They can be defined on Kleisli categories of commutative monads.
\end{abstract}


%\tableofcontents

\section{Introduction}

When dealing with a specific kind of categories one of the first questions we ask is `What kind of limits and colimits they have?'. Both operations, if they exist, are defined via universal properties, thus they are unique up to an isomorphism. However, in many circumstances we have interesting operations on a category which are not given by universal properties, yet in a given context they might be very useful. For example, one can equip a category with tensor product making it into a monoidal category \cite{B}, \cite{ML}, \cite{CWM}. Such a monoidal structure, if it exists, does not need to be unique in any sense. If we deal with categories of models of first order theories, we can equip them with ultraproduct operations \cite{L}.  These ultraproduct operations, even if they are not given by universal properties, have received a very fruitful categorical treatment (c.f. \cite{MM1}, \cite{MM2}) and proved to be useful in definability theory of first order logic, see \cite{Z}, \cite{MM3}.

The aim of this paper is to give a categorical treatment of P\l onka sums (c.f. \cite{Pl}). Originally, P\l onka sums were defined as operations on categories of algebras of regular equational theories with arites being semilattices. They are related to the (strong) sup-semilattice decomposition of semigroups \cite{T}, \cite{JLM}. Although the later perspective puts more emphasis on decomposing of a given algebra into simpler pieces rather than building more complicated algebras from simpler ones. If $L$ is a sup-semilattice considered as a (posetal) category and $\cR$ is a monad corresponding to a regular equational theory, i.e. a semi-analytic monad (c.f. \cite{SZ2}), then $L$-indexed P\l onka sum is a functor $\bigsqcup_L: \cEM(R)^L\lra \cEM(R)$, i.e. an operation on the Eilenberg-Moore category of the monad $\cR$. For $F:L\ra \cEM(\cR)$ the operation associates an $\cR$-algebra whose universe is the coproduct $\coprod_{l\in L} F(l)$ (in $Set$) and we define the operations on this universe using the transition homomorphisms $F(l\leq l'):F(l)\ra F(l')$ between the $\cR$-algebras. To calculate the value of an $n$-ary operation $f$ on elements $a_i\in F(l_i)$ for $i = 1,\ldots, n$ we first move those elements $a_i$ to a common place, that is to the algebra $F(\bigvee_i l_i)$ applying $\cR$-homomorphisms $F(l_i\leq \bigvee_i l_i)$ to elements $a_i$ and then we apply the operation $f$ to those moved elements inside the $\cR$-algebra $F(\bigvee_i l_i)$.

In categorical terms P\l onka sums are functors induced by lax morphisms of monads $(\bigsqcup_C,\phi): \hat{\cT}\ra \cT$ whose functor part is coproduct operation
\[ \bigsqcup_C: Set^C \lra Set\]
so that
\[ F \mapsto \coprod_{c\in C} F(c) \]
where $\hat{\cT}$ is the lift of the monad $\cT$ to the functor category $Set^C$ and $C$ is a small category which is the arity of the operation. The problem is what kind of monads $\cT$ and what kind of categories $C$ we should consider to get such lax morphism of monads. We show that in case of semi-analytic monads the natural choice for these arities are the categories of regular polynomials over sup-semilattices, i.e. the algebras for the terminal semi-analytic monads.  In case of analytic monads the natural choice for these arities are the categories of linear polynomials over commutative monoids, i.e. the algebras for the terminal analytic monads. There are also natural infinitary generalizations of these. In each case, the category of the monads in question (analytic, semi-analytic, and their generalizations) is a coreflexive subcategory in the category of all monads on $Set$ and the P\l onka sums have arities being categories of some kind of polynomials over the algebras for the terminal monad in this subcategory. Such P\l onka sums can be considered as an additional structure on a category. We show that the preservation of P\l onka sums by a functor between categories of algebras ensures that the morphism of monads that induced it belonged to the appropriate subcategory of monads (c.f. Theorems \ref{thm-characterization-semi-cart}, \ref{thm-characterization-semi-cart-infty}, \ref{thm-characterization-weakly-cart}).

The reason why P\l onka sums work well for semi-analytic monads is that in the corresponding theories there is a good notion of {\em occurrence of a variable in a term}. By this we mean that if two terms are equivalent modulo such a theory the same variables occur in both of them. For analytic monads P\l onka sums work, in a sense, even better (i.e. the arities can be categories of linear polynomials over monoids). This is due to the fact that in the corresponding equational theories there is a good notion of {\em number of occurrences of a variable in a term}. This means that if two terms are equivalent modulo such a theory, each variable occurs in each of them the same number of times.

As we said, the terminal objects in categories of analytic and semi-analytic monads play an important role. They have yet another property with respect to other monads in the respective categories. The commutative monoid monad distributes over any analytic monad in a canonical way and the sup-semilattice monad distributes over any semi-analytic monad in a canonical way, as well.

Looking at P\l onka sums from that abstract point of view, it is clear that there can be, at least in principle, dual operations of P\l onka products on the categories of Kleisli algebras for monads. Such monads seem to be even more rare, but surprisingly, this kind of phenomena did already appear in the literature (although dressed in a different setting). We will show that for any commutative monad $\cT$ there are operations of P\l onka product with arities being finite sets. Note that the category of finitary commutative monad is a reflective subcategory of the category of all finitary monads.

The paper is organized as follows. In Section \ref{sec-corefl-cats} we describe some subcategories of endofunctors on $Set$ and categories of monads over such categories. The only new thing in this section is the characterization of the infinitary generalizations of the category of semi-analytic functors and monads.
In Section \ref{sec-setup-for-plonka} we describe the general setup for both P\l onka sums and products. The categories of regular and linear polynomials over algebras for semi-analytic and analytic monads are described in Sections \ref{sec-reg-poly} and \ref{sec-lin-poly}, respectively. P\l onka sums on the categories of algebras for semi-analytic monads and their infinitary generalizations are described in Section \ref{sec-reg-plonka_sum}. P\l onka sums on the categories of algebras for analytic monads and their infinitary generalizations are described in Section \ref{sec-an-plonka-sums}. In Section \ref{sec-examples} some examples of concrete P\l onka sums are presented. In Sections \ref{sec-plonka-prod}  P\l onka products are discussed shortly. Finally, in Section \ref{sec-dist-laws} we discuss some distributive laws and the properties of the composed monads.

\subsection*{Preliminary notions and notation}\label{sec-prelim-notions}

In the paper we shall use category theory as well as 2-category theory language. We use the theory of monads in abstract setting of $2$-categories, including Eilenberg-Moore and Kleisli objects (c.f. \cite{St}), theory of monoidal categories and monoidal monads (c.f. \cite{CWM}).

$\o$ denotes the set of natural numbers. A cardinal number is the least ordinal of the given cardinality. Thus any $n\in\o$ is a cardinal number. If $\alpha$ is a cardinal number, then $(\alpha]=\{ 1, \ldots, \alpha\}$. $S_\alpha$ is the symmetric group, i.e. the group of permutations of the set $(\alpha]$. We write $X^\alpha$ for the set of functions from $(\alpha]$ to $X$.

\section{Coreflective subcategories of the categories of monads}\label{sec-corefl-cats}

The categories of (finitary) analytic monads $\AMnd$ and semi-analytic monads $\SanMnd$ are subcategories of the category $\Mnd$  of finitary monads on $Set$. They are categories of monoids in monoidal categories $\An$, $\San$, $\End$ of analytic, semi-analytic, and all finitary endo-functors on $Set$, respectively. The analytic functors were introduced in \cite{J}, the semi-analytic monads where introduced in \cite{M} under the name of collection monads. The categories $\An$, $\San$, $\AMnd$ and $\SanMnd$ where studied extensively in \cite{SZ1}, \cite{SZ2}.

Both kinds of functors have two characterizations: one abstract and one very concrete. Let $\S$ be the skeleton of the category of finite sets and surjection whose objects are sets $\{ 1, \ldots, n\}$, for $n\in \o$, $\B$ its subcategory with the same objects whose morphisms are bijections only.

An analytic functor in $\An$ is finitary endofunctor on $Set$ weakly preserving wide pullbacks. More concretely, $\cA$ is an analytic functor iff there is a functor $A:\B\ra Set$ and a natural isomorphism
\[   \cA(X) \cong \sum_{n\in\o} X^n\otimes_n A_n \]
where $X$ is a set and $X^n\otimes_n A_n$ is the quotient of the product of a right $S_n$-set $X^n$ with the left $S_n$-set $A_n$.

A semi-analytic functor  $\San$ is finitary endofunctor on $Set$ preserving pullbacks along monos. More concretely, $\cR$ is a semi-analytic functor iff there is a functor $R:\S\ra Set$ and a natural isomorphism
\[   \cR(X) \cong \sum_{n\in\o} X^n\otimes_n R_n \]
where $X$ is a set and $\left[\begin{array}{c} X \\ n  \end{array}\right]\otimes_n R_n$ is the quotient of the product of a right $S_n$-set $\left[\begin{array}{c} X \\ n  \end{array}\right]$ of injections from $(n]$ to $X$ with the left $S_n$-set $R_n$.  For more details see \cite{SZ2}.
The categories $\An$ and $\San$ are coreflective in $\Mnd$. If $\cT$ is a finitary endofunctor in $\End$ then its reflection in $\An$ is, for any set $X$, given by
\[ an(\cT)(X) = \sum_{n\in \o} X^n\otimes_n \cT(n) \]
 and in $\San$ is given by
\[ san(\cT)(X) = \sum_{n\in \o} \left[\begin{array}{c} X \\ n  \end{array}\right]\otimes_n \cT(n) \]
The categories $\An$ and $\San$ have terminal objects, the functors
\[ \sum_{n\in \o} X^n\otimes_n 1 \hskip 5mm {\rm and} \hskip 5mm \sum_{n\in \o} \left[\begin{array}{c} X \\ n  \end{array}\right]\otimes_n 1 \]
respectively. The morphisms in $\An$ are weakly cartesian natural transformations and the morphisms in $\San$ are semi-cartesian natural transformations, i.e. the natural transformations such that the commuting naturality squares for injections are pullbacks.

The categories $\An$, $\San$, $\End$ are (strictly) monoidal and the inclusion functors $\An\ra\San \ra \End$ are strictly monoidal functors. Thus the right adjoints $an$ and $san$ are lax monoidal functors. In particular, they preserve terminal object and so does the induced right adjoint functors between the categories of monoids, i.e. suitable monads. The terminal monad in $\AMnd$ is the monad for commutative monoids $\cC$, and the terminal monad in $\SanMnd$ is the finite power-set monad i.e. the monad for sup-semilattices $\cL.$

These concepts have various infinitary generalizations. We can replace the category of finitary endofunctors $\End$ by the category $\End_\kappa$ of $\kappa$-accessible endofunctor, where $\kappa$ is an infinite regular cardinal or $\infty$ ($\End_\infty$ is the category of all accessible functors), or even the category $\MND$ of all endofunctors on $Set$. For each category $\End_\kappa$, $\End_\infty$, $\END$ there is a corresponding (concrete) notions of analytic and semi-analytic functor ($\An_\kappa$, $\An_\infty$, $\AN$, $\San_\kappa$, $\San_\infty$, $\SAN$) and monad ($\AMnd_\kappa$, $\AMnd_\infty$, $\ANMND$, $\SanMnd_\kappa$, $\SanMnd_\infty$, $\SANMND$).

 Let$\B_\kappa$ denote the category of cardinal numbers smaller than $\kappa$ with bijections. An object $\cA$ of $\An_\kappa$ is, up to an isomorphism, given by a functor $A: \B_\kappa\ra Set$ so that for any set $X$
\[ \cA(X) \cong \sum_{\alpha\in \Card_\kappa} X^\alpha\otimes_\alpha A_\alpha \]
where $\Card_\kappa$ is the set of objects of $\B_\kappa$, i.e. the set of all cardinal numbers smaller than $\kappa$. The morphisms in $\An_\kappa$ are those induced as above by morphisms in $Set^{\B_\kappa}$, i.e. the categories $\An_\kappa$ and $Set^{\B_\kappa}$ with $\kappa$ regular are equivalent.
It can be shown that the objects of $\An_\kappa$ are $\kappa$-accessible functors that preserve wide pullbacks (c.f. \cite{AV}), and that the morphisms in this category are weakly cartesian natural transformations. The categories $\An_\infty$ and $\AN$ are the same and it is the 'sum' of all the categories $\An_\kappa$ with $\kappa$ regular cardinal.

The infinitary generalizations of semi-analytic functors haven't been considered yet, but they are also quite natural. Let $\S_\kappa$ denote the category of cardinal numbers smaller than $\kappa$ with surjections, $\S_\infty$ the category of all cardinal numbers with surjections. An object $\cR$ of $\San_\kappa$ is, up to an isomorphism, given by a functor $R: \S_\kappa\ra Set$ so that for any set $X$
\[ \cR(X) \cong \sum_{\alpha\in \Card_\kappa} \left[\begin{array}{c} X \\ \alpha  \end{array}\right]\otimes_\alpha R_\alpha \]
The morphisms in $\San_\kappa$ are induced by morphisms in $Set^{\S_\kappa}$ so that the categories $\San_\kappa$ and $Set^{\S_\kappa}$ with $\kappa$ regular cardinal are equivalent.

The category $\San_\infty$ is the 'sum' of all the categories $\San_\kappa$ with $\kappa$ regular cardinal.
An object $\cR$ of $\SAN$ is, up to an isomorphism, given by a functor $R: \S_\infty\ra Set$ so that for any set $X$
\begin{equation} \label{semi-an-fu}
\cR(X) \cong \sum_{\alpha\in \Card} \left[\begin{array}{c} X \\ \alpha  \end{array}\right]\otimes_\alpha R_\alpha
\end{equation}
$\Card$ is the class of all cardinal numbers\footnote{By a cardinal number we mean the least ordinal number of the given cardinality.}. Note that, this sum is well defined as for $\alpha > card(X)$ the set $ \left[\begin{array}{c} X \\ \alpha  \end{array}\right]$ of injections from $(\alpha]$ to $X$ is empty. The categories $\San_\infty$ and $Set^{\S_\infty}$ are equivalent. The following theorem provides a characterization of the above categories.

\begin{theorem} \label{semi-analytic-char}
The objects of  the category $\San_\kappa$, where $\kappa$ is either regular cardinal or $\infty$, are $\kappa$-accessible functors that preserve pullbacks along monos and wide pullbacks of monos. The objects of  the category $\SAN$ are functors that preserve pullbacks along monos and wide pullbacks of monos. The morphisms in these categories are semi-cartesian natural transformations.
\end{theorem}

{\em Proof.} We shall prove the characterization for $\San_\kappa$ where $\kappa$ is any regular cardinal. The other two cases are so similar that can be easily deduced.

Thus, by an argument analogous to the one given in the proof of Proposition 2.1 of \cite{SZ2}, we need to identify the image of Kan extension functor $Set^{\S_\kappa}\ra \END$. This is a `lift' of the proof of Theorem 2.2 of \cite{SZ2} from $\kappa=\o$ to any regular cardinal $\kappa$. On one hand, if an endofunctor $\cR$ on $Set$ is defined by the formula (\ref{semi-an-fu}), then it easy to see that it preserves wide pullbacks of monos. On the other hand, if an endofunctor $F$ on $Set$ is $\kappa$-accessible and preserves wide pullbacks of monos then for any set $X$ and any $x\in F(X)$ there is the least cardinal $\alpha$ necessarily smaller then $\kappa$ and the least subobject $F:(\alpha]\lra X$ and $y\in F(\alpha]$ such that $F(f)(y)=x$. With these observations the rest can be easily deduced form the proof of Theorem 2.2. in \cite{SZ2}. $\boxempty$
\vskip 2mm

 Having the above characterization it is clear now that all the categories of functors considered above are closed under composition of their objects. For the explicite formulas for the coefficient functors of the composed functors in categories $\An$ and $\San$ the reader can consult \cite{SZ1} and \cite{SZ2}. Thus we have described a diagram of (strictly monoidal) categories and functors.
\begin{center} \xext=2100 \yext=1150
\begin{picture}(\xext,\yext)(\xoff,\yoff)
%lower squares
\setsqparms[1`-1`-1`1;700`500]
 \putsquare(0,50)[\San`\San_\kappa`\An`\An_\kappa;```]
\setsqparms[1`0`-1`1;700`500]
 \putsquare(700,50)[\phantom{\San_\kappa}`\San_\infty`\phantom{\An_\kappa}`\An_\infty;```]
  \putsquare(1400,50)[\phantom{\San_\infty}`\SAN`\phantom{\An_\infty}`\AN;```=]
%upper squares
 \setsqparms[1`-1`-1`0;700`500]
 \putsquare(0,550)[\End`\End_\kappa`\phantom{\San}`\phantom{\San_\kappa};```]
\setsqparms[1`0`-1`0;700`500]
 \putsquare(700,550)[\phantom{\End_\kappa}`\End_\infty`\phantom{\San_\kappa}`\phantom{\San_\infty};```]
  \putsquare(1400,550)[\phantom{\End_\infty}`\END`\phantom{\San_\infty}`\phantom{\SAN};```]
 \end{picture}
\end{center}
All the functors are (strictly monoidal) inclusions. Horizontal inclusions are full. All the categories but $\An_\infty = \AN$ and $\San_\infty$ are coreflective in all the categories that contain them. All the categories but $\An_\infty = \AN$ and $\San_\infty$ have terminal objects which are the values of the right adjoints to the inclusion on the terminal object in $\END$.

The above diagram lifts to the diagram of monoids in those categories and we get the following categories of monads with functors being again inclusions
\begin{center} \xext=2100 \yext=1150
\begin{picture}(\xext,\yext)(\xoff,\yoff)
%lower squares
\setsqparms[1`-1`-1`1;700`500]
 \putsquare(0,50)[\SanMnd`\SanMnd_\kappa`\AMnd`\AMnd_\kappa;```]
\setsqparms[1`0`-1`1;700`500]
 \putsquare(700,50)[\phantom{\SanMnd_\kappa}`\SanMnd_\infty`\phantom{\AMnd_\kappa}`\AMnd_\infty;```]
  \putsquare(1400,50)[\phantom{\SanMnd_\infty}`\SANMND`\phantom{\AMnd_\infty}`\ANMND;```=]
%upper squares
 \setsqparms[1`-1`-1`0;700`500]
 \putsquare(0,550)[\Mnd`\Mnd_\kappa`\phantom{\SanMnd}`\phantom{\SanMnd_\kappa};```]
\setsqparms[1`0`-1`0;700`500]
 \putsquare(700,550)[\phantom{\Mnd_\kappa}`\Mnd_\infty`\phantom{\SanMnd_\kappa}`\phantom{\SanMnd_\infty};```]
  \putsquare(1400,550)[\phantom{\Mnd_\infty}`\MND`\phantom{\SanMnd_\infty}`\phantom{\SANMND};```]
 \end{picture}
\end{center}
If a monoidal category has a terminal object $1$, then there is a unique structure of a monoid on $1$ and this monoid is the terminal monoid. Therefore all the above categories of monads but $\AMnd_\infty = \ANMND$ and $\SanMnd_\infty$  have terminal objects. The terminal object in $\AMnd$ is the monad for commutative monoids and the terminal object in $\AMnd_\kappa$ is its 'less than $\kappa$' version. The terminal object in $\SanMnd$ is the monad for sup-semilattices and the terminal object in $\SanMnd_\kappa$ is its 'less than $\kappa$' version. The terminal object in $\SANMND$ is the power-set monad i.e. the monad for suplattices.

\section{A general setup for P\l onka sums and products}\label{sec-setup-for-plonka}

Let $\CAT$ denote the $2$-category of possibly large categories, functors, and natural transformations. Let $\MND_{lax}$ and $\MND_{oplax}$ denote the $2$-categories of monads in $\CAT$ with lax and oplax morphism of monads, respectively, and with suitable notion of transformation of such morphism of monads.
Let $\cEM : \MND_{lax} \lra \CAT$ denote the $2$-functor of Eilenberg-Moore object, and  $\cK : \MND_{oplax} \lra \CAT$ denote the $2$-functor of Kleisli object (c.f. \cite{St}).

Let $\cS=(\cS,\eta,\mu)$ be a monad on a complete and cocomplete category $\cA$. If $C$ is a small category than $\cS$ lifts to a monad $\hat{\cS}$ on the functor category $\cA^C$ so that the diagonal functor
\[ \delta_\cA: \cA \lra \cA^C \]
is a strict morphism of monads $\cS\ra \hat{\cS}$. As $\delta_\cA$ has both adjoints $Colim_{\cA,C}\dashv \delta_\cA\dashv Lim_{\cA,C}$ it induces an oplax morphism of monads $(Colim_{\cA,C}, \psi_c) :  \hat{\cS}\ra\cS$ and lax morphism of monads $(Lim_{\cA,C}, \psi_l) :  \hat{\cS}\ra\cS$.

The exponentiation 2-functor $(-)^C : \CAT\ra\CAT$ commutes with Eilenberg-Moore objects in $\CAT$. Thus the functor $\cEM(\delta_\cA)$ between Eilenberg-Moore categories induced by the strict morphism of monads $\delta_\cA$ factorizes via the diagonal functor $\delta_{\cEM(\cS)}$ as follows
\begin{center} \xext=1400 \yext=100
\begin{picture}(\xext,\yext)(\xoff,\yoff)
  \putmorphism(0,50)(1,0)[\cEM(\cS)`\cEM(\cS)^C`\delta_{\cEM(\cS)}]{800}{1}a
  \putmorphism(800,50)(1,0)[\phantom{\cEM(\cS)^C}`\cEM(\hat{\cS})`\simeq]{600}{1}a
\end{picture}
\end{center}
As $(Lim_C, \psi_l)$ is the right adjoint to $\delta_\cA$ in $\MND_{lax}$, the functor $\cEM(Lim_C, \psi_l)$ is the right adjoint to $\cEM(\delta_\cA)$ in $\CAT$ and it factorizes as
\begin{center} \xext=1400 \yext=150
\begin{picture}(\xext,\yext)(\xoff,\yoff)
  \putmorphism(0,50)(1,0)[\cEM(\hat{\cS})`\cEM(\cS)^C`\simeq]{600}{1}a
  \putmorphism(600,50)(1,0)[\phantom{\cEM(\cS)^C}`\cEM(\cS)`\overline{Lim}_C]{800}{1}a
\end{picture}
\end{center}
$\overline{Lim}_C$ is nothing but the functor of taking limits of $\cS$-algebras indexed by the small category $C$. This is a way we can explain why categories of Eilenberg-Moore algebras are complete whenever the categories over which they are defined are.

As we know, the situation with Kleisli algebras is very different. Typically, the Kleisli category is not cocomplete even if the category over which it is defined is. This is due to the fact that the exponentiation $2$-functor $(-)^C : \CAT\ra\CAT$ does not commutes with Kleisli objects, in general. It does if the index category $C$ is discrete. So the above argument can be repeated only in case $C$ is a discrete category, i.e. in such case the functor $\cK(Colim_{\cA,C}, \psi_c)$, the left adjoint to $\cK(\delta_\cA)$, factorizes as
\begin{center} \xext=1400 \yext=150
\begin{picture}(\xext,\yext)(\xoff,\yoff)
  \putmorphism(0,50)(1,0)[\cK(\hat{\cS})`\cK(\cS)^C`\simeq]{600}{1}a
  \putmorphism(600,50)(1,0)[\phantom{\cK(\cS)^C}`\cK(\cS)`\overline{Colim}_C]{800}{1}a
\end{picture}
\end{center}
and $\overline{Colim}_C$ is the usual coproduct of free algebras. This is a way we can explain why categories of Kleisli algebras have coproducts, but not all colimits in general, whenever the categories over which they are defined have either coproducts or even all colimits.

The natural transformation $\psi_c$ is a natural isomorphism i.e. $(Colim_C, \psi_c)$ is a strong morphism of monads iff $\cS$ preserves $C$ indexed colimits. If $\cA=Set$, $\cS$ is finitary and this holds even for (finite) discrete categories $C$, it implies that all the operations in the equational theory corresponding to the monad $\cS$ are unary. Thus such monads are very rare and mostly uninteresting. However, if a monad $\cS$ has some additional good properties (like being analytic or semi-analytic on $\cA=Set$) it may happen that for some small category $C$ the functor $\bigsqcup_C : \cA^C \lra  \cA$, the composition of functors
\[ \cA^C \stackrel{i^*}{\lra}  \cA^{ob(C)} \stackrel{\coprod}{\lra}  \cA \]
induced by functors between small categories
$C\stackrel{i}{\lla} ob(C) \stackrel{!}{\lra} 1$, can be equipped with a natural transformation $\lambda_C : \cS\circ \bigsqcup_C \lra  \bigsqcup_C \circ \hat{\cS}$ so that  $(\bigsqcup_C,\lambda_C) : \hat{\cS}\ra \cS$ is a lax morphism of monads. Each such lax morphism of monads induces in turn an operation on the category of Eilenberg-Moore algebras
\[ \bigsqcup{\!_C} : \cEM(\cS)^C\cong \cEM(\hat{\cS})\lra \cEM(\cS) \]
In a special case these operations are what is called P\l onka sum on the category of algebras of a regular equational theory. If $\cS$ is a semigroup monad, $C$ is a sup-semilattice and $\bF: C \ra \cEM(S)$ is a functor, then $\bF$ is what is called a {\em strong sup-semilattice decomposition of the semigroup} $\bigsqcup_C(\bF)$. Thus such lax morphisms of monads induce additional operations on the categories of Eilenberg-Moore algebras that we shall call (generalized) P\l onka sums.

Dually, the natural transformation $\psi_l$ is an isomorphism iff $\cS$ preserves $C$ indexed limits. Such monads are even more rare. But again, if a monad $\cS$ has some additional good properties, it may happen that for some small categories $C$ the functor $\bigsqcap_C : \cA^C \lra  \cA$, the composition of functors
\[ \cA^C \stackrel{i^*}{\lra}  \cA^{ob(C)} \stackrel{\prod}{\lra}  \cA \]
induced by functors between small categories
$C\stackrel{i}{\lla} ob(C) \stackrel{!}{\lra} 1$, can be equipped with a natural transformation $\rho_C :  \bigsqcap_C \circ \hat{\cS}\lra \cS\circ \bigsqcap_C$ so that  $(\bigsqcap_C,\rho_C) : \hat{\cS}\ra \cS$ is a oplax morphism of monads. Each such oplax morphism of monads induces in turn an operation on the category of Kleisli algebras
\[ \bigsqcap{\!_C} : \cK(\cS)^C\cong \cK(\hat{\cS})\lra \cK(\cS) \]
By analogy we shall call such operations P\l onka products.

\vskip 2mm
{\em Remark.} For any ultrafilter $U$ on any set $I$ we have an ultraproduct functor $Set^I\stackrel{[U]}{\lra} Set$, see \cite{MM1} for details. If $\cT$ is a finitary monad on $Set$ and $\hat{\cS}$ its lift to $Set^I$ then $[U]:\hat{\cS}\lra \cS$ is a strong monad morphism, as a consequence of \L o\'{s} theorem. It induces the ultraproduct operation on the category of $\cS$-algebras.

\section{Category of regular polynomials over an algebra}\label{sec-reg-poly}

Let $\cR=(\cR,\eta,\mu)$ be a semi-analytic monad (c.f. \cite{SZ2}), $R:\S \ra Set$ a functor such that, for any set $X$
\[ \cR(X)= \sum_{n\in\o} \left[\begin{array}{c} X \\ n  \end{array}\right] \otimes_n R_n \]
By \cite{SZ2}, $R$ is the functor part of a regular operad, i.e. a monoid in the monoidal category $Set^\S$ with the substitution tensor.  We define a functor
\[ CP_r^\cR : \cEM(\cR) \lra Cat \]
associating to an $\cR$-algebra $(A,\alpha:\cR(A)\ra A)$ a category of regular polynomials $CP_r^\cR(A,\alpha)$ as follows. The objects of $CP_r^\cR(A,\alpha)$  are elements of $A$. A morphism in $CP_r^\cR(A,\alpha)$ is an equivalence class of triples
\[ [\vec{a}, i, r ]_\sim : \vec{a}(i) \ra \alpha([\vec{a},r]_\sim) \]
where $\vec{a}:(n]\ra A$ is an injection, $i\in(n]$, $r\in R_n$, for some $n\in \o$. Note that $[\vec{a},r]_\sim$ is an element of $\cR(A)$. We identify triples
 \[ \lk \vec{a}\circ\sigma, i, r \rk \sim \lk \vec{a}, \sigma(i), R(\sigma)(r) \rk \]
where $\sigma\in S_n$. The identity morphism is
\[  [a, 1, \iota ]_\sim : a \ra \alpha([a,\iota]_\sim)=a \]
where $\iota\in R_1$ is the unit of the regular operad $R$. The composition is defined by the substitution of regular terms into regular terms possibly with normalization. In detail, it is defined as follows
\begin{center}
\xext=2400 \yext=280
\begin{picture}(\xext,\yext)(\xoff,\yoff)
\putmorphism(0,50)(1,0)[\vec{a}(i)`{\alpha(\vec{a},r)=\vec{a}'(j)}`[\vec{a}, i, r {]}_\sim]{1200}{1}b
\putmorphism(1200,50)(1,0)[\phantom{\alpha(\vec{a},r)=\vec{a}'(j)}`\alpha(\vec{a}',r')`[\vec{a}', j, r' {]}_\sim]{1200}{1}b
\putmorphism(0,180)(1,0)[\phantom{\vec{a}(i)}`\phantom{\alpha(\vec{a}',r)}`[a'', i'', r''  {]}_\sim]{2400}{1}a
\end{picture}
\end{center}
where $r\in R_n$, $r'\in R_m$; the function $\vec{a}'':(k]\ra A$ is the injection part for surjection-injection factorization $\vec{a}''\circ s:(n+m-1] \ra A$
of the function    $\vec{a}'(j\backslash\vec{a}):(n+m-1] \ra A$ such that
\[ \vec{a}'(j\backslash\vec{a})(l)\;\;= \;\; \left\{ \begin{array}{ll}
                \vec{a}'(l)   & \mbox{ if } 1\leq l<j \\
                \vec{a}(l-j+1)   & \mbox{ if } j\leq l <n+j \\
                \vec{a}'(l-n+1) & \mbox{ if } n+j\leq l < n+m .
                                    \end{array}
                            \right. \]
Thus $\vec{a}'(j\backslash\vec{a})$ replaces $j$ in the domain of $\vec{a}'$ by the whole function $\vec{a}$. Such function might not be an injection and $\vec{a}''$ is the injection part of it. $i''=s(i+j-1)\in (k]$, i.e. $i''$ is the element in $(k]$ that correspond to $i\in (n]$. $r''=R(s)((\iota,\ldots,r,\ldots,\iota)\ast r')$, i.e. the value under the action of $s:(n+m-]\ra (k]$ on the composition in the regular operad $R$ of $r'$ with $r$ placed into the $j$'s place. This ends the definition of the category $CP_r^\cR(A,\alpha)$.

A homomorphism $h:(A,\alpha)\lra (A',\alpha')$  induces a functor \[ CP_r^\cR(h) :CP_r^\cR(A,\alpha)\lra CP_r^\cR(A',\alpha')\]
so that the morphism $[\vec{a}, i, r ]_\sim : \vec{a}(i) \ra \alpha(\vec{a},r)$ is sent to
\[ [\vec{a}', s(i), R(s)(r) ]_\sim : \vec{a}'(s(i))=h(\vec{a}(i)) \ra h(\alpha(\vec{a},r))= \alpha'(\vec{a}',R(s)(r)) \]
where $\vec{a}'\circ s$ is surjection-injection factorization of $h\circ \vec{a}$.

We note for the record

\begin{fact}\label{reg-polynomials} For any semi-analytic monad $\cR$,
the functor
\[ CP_r^\cR : \cEM(\cR) \lra Cat \]
associating to $\cR$-algebras their categories of regular polynomials is well defined. $\boxempty$
\end{fact}

\vskip 2mm
{\em Remark.} One can describe the category of polynomials $CP_r^\cR(A,\alpha)$ in terms of the regular equational theory $T_\cR$ corresponding to the monad $\cR$ as follows. Its objects are elements of the algebra $(A,\alpha)$.
A morphism from $a$ to $b$ is given by
\begin{enumerate}
  \item a regular term of the theory $T_\cR$ in $n$ variable $r(x_1,\ldots,x_n)$ (all the variables necessarily explicitly occur in $r$);
  \item a certain variable number $i\in(n]$;
  \item an injective function $\vec{a}:(n]\ra A$ interpreting variable occurring in $r(x_1,\ldots,x_n)$;
  \item the domain of the morphism $a$ is equal to the interpretation of the $i$'th variable $\vec{a}(i)$;
  \item the codomain of the morphism $b$ is equal to the value of the term $r(x_1,\ldots,x_n)$ under the interpretation of the variables $\vec{a}$.
\end{enumerate}
The composition of such morphisms is defined as (the most reasonable) substitution.

The reason why this definition works for semi-analytic monads is, as we explained in the introduction, that in the corresponding regular equational theories there is a good notion of occurrence of a variable in a term.
\vskip 2mm

If $\tau: \cR\ra \cV$ is a morphism of semi-analytic monads, then the functors associating categories of regular polynomials to those monads are related as follows. The morphism $\tau$ induces a 'forgetful' functor $\cEM(\tau): \cEM(\cV)\ra \cEM(\cR)$ and we have a natural transformation
\[  \gamma^\tau:  CP^\cR_r\circ \cEM(\tau) \lra CP^\cV_r : \cEM(\cV)\lra Cat\]
defined for a $\cV$-algebra $(B,\beta)$ a functor
\[ \gamma^\tau_{(B,\beta)}: CP^\cR_r(B,\beta\circ\tau) \lra CP^\cV_r(B,\beta) \]
constant on object, and sending morphism $[\vec{a}, i, r ]_\sim :\vec{a}(i)\lra \beta\circ\tau(\vec{a},r)$ to the morphism
\[ [\vec{a}, i, \tau_n(r) ]_\sim : \vec{a}(i)\lra \beta(\vec{a},\tau_n(v))  \]
where $n\in\o$, $r\in R_n$ and $\tau_n: R_n\ra V_n$.

\subsection*{Infinitary case}

It should be clear how to define the functor of associating to algebras for a monad $\cR$ in $\SanMnd_\kappa$, $\SanMnd_\infty$ or $\SANMND$ the category of (suitable infinitary) regular polynomials
\[ CP_{\infty r}^\cR : \cEM \lra Cat \]
in analogy with the finitary case. Note that we never consider regular polynomials of arity that exceeds the cardinality of the algebra. The details are left for the reader.

%\newpage
\section{P\l onka sums of algebras for a semi-analytic monad}\label{sec-reg-plonka_sum}

In this section we shall study operations of generalized P\l onka sums on the categories of algebras for semi-analytic monad indexed by categories of regular polynomials of algebras for another semi-analytic monads.

Let $\pi : \cR\ra \cT$ be a morphism of semi-analytic monads, defined by a natural transformation (denoted by the same letter) $\tau:R\ra T$ in $Set^\S$. Let $(A,\alpha)$ be a $\cT$-algebra. Let us denote  the category of regular polynomials ${CP_r(A,\alpha)}$ as $\bA$, for short. The monad $\hat{\cR}$ is the lift of the monad $\cR$ to the category $Set^{\bA}$. We shall define a lax morphism of monads induced by $\pi$ and $\cT$-algebra $(A,\alpha)$
\[ (\bigsqcup{\!_{(A,\alpha)}},\lambda^{\pi,(A,\alpha)}) :  \hat{\cR} \ra \cR \]
We usually drop superscripts  $(A,\alpha)$ and write $\cR$ instead of $\pi$ in $\lambda^{\pi,(A,\alpha)}$ when it does not lead to a confusion. We also write $\bigsqcup{\!_{A}}$ rather than $\bigsqcup{\!_{(A,\alpha)}}$. Let $\bF:\bA\ra \cEM(\cR)$ be a functor, and $F:\bA\ra Set$ the composition of $\bF$ with the forgetful functor. We shall define the component
\[ \lambda^\cR_F: \cR(\coprod_{a\in A}F(a)) \lra \coprod_{a\in A}\cR(F(a)) \]
of the natural transformation $\lambda^\cR$. Let $[\vec{x},r]_\sim \in \cR(\coprod_{a\in A}F(a))$, where $\vec{x}:(n]\ra \coprod_{a\in A} F(a)$ is an injection, $r\in R_n$. Let $p:\coprod_{a\in A} F(a)\ra A$ be the projection from the coproduct to the index set. Let $\vec{a}\circ s$ be a surjection-injection factorization of $p\circ \vec{x}$ as in the diagram below. Moreover, we have a unique function $\vec{x}'$ making the triangle in the middle commute, where $\kappa_a: F(a)\ra \coprod_{a\in A} F(a)$ is the injection into coproduct. To explain the right hand square in the diagram below, note that $[\vec{a},R(s)(r)]\in \cR(A)$ and hence
\[ \pi([\vec{a},R(s)(r)])=[\vec{a},\pi_m(R(s)(r))]\in \cT(A)\]
We put
\[ b=\alpha([\vec{a},\pi_m(R(s)(r))]) \]
and, for $i\in (n]$, we have a morphism
\[ \psi_i = [\vec{a},s(i),\pi_m(R(s)(r))]: \vec{a}(s(i))\lra b\]
in the category $\bA$. The morphism $\vec{z}\circ g$ is a surjection-injection factorization of $[F(\psi_i)]_{i\in(n]}\circ\vec{x}'$. This explains the construction of the diagram
\begin{center} \xext=3000 \yext=700
\begin{picture}(\xext,\yext)(\xoff,\yoff)
\setsqparms[0`-1`0`-1;1200`500]
\putsquare(0,50)[A``{(m]}`{(n]};`\vec{a}``s]

\setsqparms[0`0`-1`1;1800`500]
\putsquare(1200,50)[`F(b)`\phantom{(n]}`{(k]};``\vec{z}`g]

\putmorphism(0,550)(1,0)[\phantom{A}`\coprod_{a\in A} F(a)`p]{600}{-1}a
\putmorphism(600,550)(1,0)[\phantom{\coprod_{a\in A} F(a)}`\coprod_{i\in (n]} F(\vec{a}(s(i)))`{[\kappa_{\vec{a}(s(i))}]}_{i\in(n]}]{1300}{-1}a
\putmorphism(1900,550)(1,0)[\phantom{\coprod_{i\in (n]} F(\vec{a}(s(i)))}`\phantom{F(b)}`[F(\psi_i){]}_{i\in(n]}]{1100}{1}a

          \put(1130,140){\vector(-1,1){320}}
          \put(1240,140){\vector(1,1){320}}

\put(900,200){${\vec{x}}$}
\put(1400,200){${\vec{x}'}$}
 \end{picture}
\end{center}
Finally, we put
\[ \lambda^\cR_F([\vec{x},r]_\sim) = [\vec{z},R(g)(r)]_\sim \]

A simple verification shows

\begin{proposition}\label{preservation_polonka}
$(\bigsqcup{\!_{A}},\lambda^\cR) :  \hat{\cR} \ra \cR$ is a lax morphism of monads. $\boxempty$
\end{proposition}

Lifting the above morphism of monads, we obtain $\cT$-indexed P\l onka sum of $\cR$-algebras, i.e. for any $\cT$-algebra $(A,\alpha)$ we obtain an operation
\begin{center} \xext=1350 \yext=750
\begin{picture}(\xext,\yext)(\xoff,\yoff)
\setsqparms[1`1`1`1;1200`550]
 \putsquare(150,50)[\cEM(\cR)^{CP_r(A,\alpha)}\cong \cEM(\hat{\cR})`\cEM(\cR)`Set^{CP_r(A,\alpha)}`Set;\bigsqcup{\!_{A}}`\cU^{CP_r(A,\alpha)}`\cU`\bigsqcup{\!_{A}}]
\end{picture}
\end{center}

The following Proposition explains how P\l onka sums interact with lax morphisms of semi-analytic monads.

\begin{proposition}\label{T-morphisms-pres-plonka} Let
\begin{center} \xext=600 \yext=320
\begin{picture}(\xext,\yext)(\xoff,\yoff)
 \settriparms[1`1`1;300]
  \putVtriangle(0,0)[\cR`\cV`\cT;\tau`\pi`\pi']
\end{picture}
\end{center}
be a commuting triangle in the category of semi-analytic monads, $(A,\alpha)$ a $\cT$-algebra. Then we have a commuting square
\begin{center} \xext=800 \yext=600
\begin{picture}(\xext,\yext)(\xoff,\yoff)
 \setsqparms[-1`1`1`-1;800`500]
 \putsquare(0,50)[\hat{\cR}`\hat{\cV}`\cR`\cV;(1,\hat{\tau})`(\bigsqcup{\!_{A}},\lambda^\cR)`(\bigsqcup{\!_{A}},\lambda^\cV)`(1,\tau)]
  \end{picture}
\end{center}
of lax morphisms of monads, where transformations $\lambda^\cR$ and $\lambda^\cV$ are induced by morphism $\pi$ and $\pi'$, respectively.
\end{proposition}

{\em Proof.} Let us fix a $\cT$-algebra $(A,\alpha)$ and $\bA=CP_r(A,\alpha)$. We shall denote the natural transformations in $Set^\S$ that give rise to $\tau$, $\pi$, $\pi'$ by the same letters.  We need to show that the square
of functors and natural transformations in $Nat(Set^\bA,Set)$
\begin{equation}\label{pres-polonka-square}\end{equation}
\begin{center} \xext=1200 \yext=500
\begin{picture}(\xext,\yext)(\xoff,\yoff)
 \setsqparms[1`-1`-1`-1;1200`500]
 \putsquare(0,50)[\bigsqcup{\!_{A}}\circ \cR`\bigsqcup{\!_{A}}\circ \cV`\cR\circ \bigsqcup{\!_{A}}`\cV\circ\bigsqcup{\!_{A}};
 \bigsqcup{\!_{A}}(\tau)`\lambda^\cR`\lambda^\cV`\tau_{(\bigsqcup{\!_{A}})}]
  \end{picture}
\end{center}
commutes.

Let $F:\bA\ra Set$ be a functor. An element $[\vec{x},r]_\sim$ of $\cR\circ \bigsqcup_A(F)$ is represented by a pair $\lk \vec{x},r\rk$ such that, for some $n\in \o$, $\vec{x}:(n]\ra \bigsqcup_{a\in A}F(a)$ is an injection and $r\in R_n$. Let $\vec{a}$, s, $\vec{x}'$, $b$, $\psi_i$, $g$, $\vec{z}$ be as in the definition of $\lambda^\cR_F([\vec{x},r]_\sim)$ above. Thus
\begin{equation}\label{eq1} \bigsqcup{\!_A}(\tau)(\lambda^\cR([\vec{x},r]_\sim)) =\tau_{F(b)}([\vec{z},R(g)(r)]_\sim)=  {[\vec{z},\tau_m(R(g)(r))]_\sim}
\end{equation}

As $\tau$ is a natural transformation and  $\pi'\circ \tau = \pi$, for $i\in (n]$, the morphism $\psi_i$ in $\bA$ is equal to the morphism
\[ \psi'_i = [\vec{a},s(i),\pi'_m(V(s)(\tau_n(r)))]_\sim: \vec{a}(s(i))\lra b \]
Thus
\begin{equation}\label{eq2} \lambda^\cV (\tau_{(\bigsqcup{\!_{A}})}( [\vec{x},r]_\sim)) = \lambda^\cV ([\vec{x},\tau_n(r)]_\sim))=
 {[\vec{z},V(s)(\tau_n(r))]_\sim}
\end{equation}
Now from (\ref{eq1}) and (\ref{eq2}) if follows that the square (\ref{pres-polonka-square}) commutes. $\boxempty$
\vskip 2mm

Lifting Proposition \ref{T-morphisms-pres-plonka} to the categories of Eilenberg-Moore algebras, with the notation as above, we obtain
that regular interpretations of regular theories induce morphisms of algebras that preserve P\l onka sums. Thus we have

\begin{corollary}
The functor between categories of algebras $\cEM(\tau) : \cEM(\cV)\ra \cEM(\cR)$ induced by the semi-analytic natural transformation $\tau:\cR\ra\cT$ preserves $\cT$-indexed P\l onka sum, i.e. for any  $\cT$-algebra $(A,\alpha)$ the square
\begin{center} \xext=1000 \yext=750
\begin{picture}(\xext,\yext)(\xoff,\yoff)
 \setsqparms[-1`1`1`-1;1000`600]
 \putsquare(0,70)[\cEM(\cR)^\bA`\cEM(\cV)^\bA`\cEM(\cR)`\cEM(\cV);\cEM(\tau)^\bA`\bigsqcup{\!_{(A,\alpha)}}`\bigsqcup{\!_{(A,\alpha)}}`\cEM(\tau)]
  \end{picture}
\end{center}
commutes, up to a canonical isomorphism, where $\bA=CP_r(A,\alpha)$.
$\boxempty$
\end{corollary}
\vskip 2mm
From Proposition \ref{T-morphisms-pres-plonka} we also have

\vskip 2mm
\begin{corollary} \label{Polnka-pres-algebras}
Let
\begin{center} \xext=300 \yext=320
\begin{picture}(\xext,\yext)(\xoff,\yoff)
 \settriparms[1`1`1;300]
  \putVtriangle(0,0)[\cR`\cV`\cT;\tau`\pi`\pi']
\end{picture}
\end{center}
be a commuting triangle in the category of semi-analytic monads, $(A,\alpha)$ a $\cT$-algebra $\bA=CP_r(A,\alpha)$. Let $\bF: \bA\ra \cEM(\cV)$ be a functor,
$F=\cU^\cV\circ \bF : \bA\ra Set$. Then the two $\cR$-algebra structures on $\coprod_{a\in A} F(a)$
\begin{center} \xext=1600 \yext=1850
\begin{picture}(\xext,\yext)(\xoff,\yoff)
 \settriparms[1`1`0;600]
  \putAtriangle(0,1220)[\cR(\coprod_aF(a))`\phantom{\coprod_a \cR(F(a))}`\phantom{\cV(\coprod_aF(a))};\lambda^\cR_F`\tau_{\coprod_aF(a)}`]

\setsqparms[0`1`1`0;1200`600]
 \putsquare(0,600)[\coprod_a \cR(F(a))`\cV(\coprod_aF(a))`\coprod_a\cV(F(a))`\coprod_a\cV(F(a));`\coprod_a\tau_{F(a)}`\lambda^\cV_F`]

 \settriparms[0`1`1;600]
 \putVtriangle(0,0)[``\coprod_aF(a);`\coprod_a\xi_a`\coprod_a\xi_a]
 \end{picture}
\end{center}
coincide.
$\boxempty$
\end{corollary}

\vskip 2mm

In order to make these sums independent of a given monad $\cT$ and particular morphisms $\pi$ and $\pi'$, we can take as $\cT$ the terminal semi-analytic monad $\cL$ of sup-semilattices. In such a way we have equipped any category of algebras for a semi-analytic monad $\cR$ with a canonical system of P\l onka sums indexed by the categories of regular polynomials over sup-semilattices. By Proposition \ref{preservation_polonka} the semi-analytic morphism of semi-analytic monads induces a functor between categories of algebras that preserves P\l onka sums. In fact, this property characterizes semi-analytic morphisms of semi-analytic monads, i.e. the converse of Corollary \ref{Polnka-pres-algebras}, and hence also Proposition \ref{preservation_polonka},  is also true. We have

\begin{theorem}\label{thm-characterization-semi-cart} Let $\tau:\cR\ra \cV$ be an arbitrary lax morphism of monads between semi-analytic monads. Then $\tau$ is a semi-cartesian iff the induced functor $\cEM(\tau)$ between categories of algebras preserves P\l onka sums.
\end{theorem}

{\em Proof.} The `only if' part is just Proposition \ref{T-morphisms-pres-plonka}.

To see the converse, for a lax morphism of monads $\tau :\cR\ra \cV$, we shall define a natural transformation $\sigma : R\ra V$ such that $\hat{\sigma}=\tau$.
Let us fix $r\in R_n$, $n\in \o$. Then $[1_{(n]},r]_\sim \in \cR(n]$ and we have
\[ \tau_{(n]}([1_{(n]},r]_\sim) = [\vec{x},v]_\sim \in \cV(n] \]
where $\vec{x}:(m]\ra (n]$ is an injection and $v\in V_m$, for some $m\in\o$. To end the proof
\begin{enumerate}
  \item we show that $m=n$ and hence $\vec{x}$ is a permutation;
  \item we put $\sigma_n(r)=V(\vec{x})(r)$ and show that $\sigma : R\ra V$ so defined is a natural transformation;
  \item finally we show that $\hat{\sigma}=\tau$.
\end{enumerate}

Let $X\subseteq (n]$ be the image of $\vec{x}$ is $(n]$. Let $\bF : \cP(n] \lra \cEM(\cV)$ be a functor sending $y\in \cP(n]$ to the one-element algebra $(\{ Y\},\xi_Y)$. Thus we have a commuting diagram
\begin{center} \xext=1100 \yext=750
\begin{picture}(\xext,\yext)(\xoff,\yoff)
 \settriparms[1`1`1;650]
  \putqtriangle(0,0)[\cP(n]`\cEM(\cV)`Set;\bF`F`\cU^\cV]
 \settriparms[1`0`1;650]
  \putptriangle(650,0)[\phantom{\cEM(\cV)}`\cEM(\cV)`\phantom{Set};\cEM(\tau)``\cU^\R]
\end{picture}
\end{center}

with $F$ being an inclusion functor. Note that $\coprod_{Y\in\cP(n]}\{ Y\} =\cP(n]$.

Let $\{ -\} : (n] \ra \cP(n]$ be the singleton morphism i.e. $\{-\}(i) = \{ i\}$, for $i\in (n]$. Since $\cEM(\tau)$ preserves P\l onka sums and$\tau$ is a natural transformation we have a commuting diagram
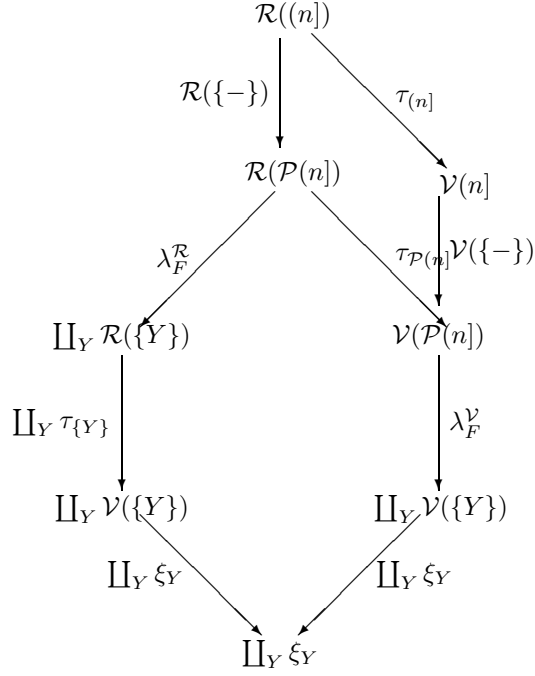
\begin{center} \xext=1200 \yext=2550
\begin{picture}(\xext,\yext)(\xoff,\yoff)
\putmorphism(600,2450)(0,-1)[\phantom{(n]}`\phantom{\cP(n]}`\cR(\{-\})]{550}{1}l
\putmorphism(1200,1850)(0,-1)[\phantom{(n]}`\phantom{\cP(n]}`\cV(\{-\})]{550}{1}r

 \settriparms[0`1`0;650]
  \putAtriangle(0,1820)[\cR((n])`\phantom{\coprod_Y \cR(\{ Y\})}`\cV(n];`\tau_{(n]}`]

 \settriparms[1`1`0;650]
  \putAtriangle(0,1220)[\cR(\cP(n])`\phantom{\coprod_Y \cR(\{ Y\})}`\phantom{\cV(\coprod_aF(a))};\lambda^\cR_F`\tau_{\cP(n]}`]

\setsqparms[0`1`1`0;1200`650]
 \putsquare(0,600)[\coprod_Y \cR(\{Y\})`\cV(\cP(n])`\coprod_Y\cV(\{ Y\})`\coprod_Y\cV(\{ Y \});`\coprod_Y\tau_{\{ Y\}}`\lambda^\cV_F`]

 \settriparms[0`1`1;600]
 \putVtriangle(0,50)[``\coprod_Y\xi_Y;`\coprod_Y\xi_Y`\coprod_Y\xi_Y]

 \end{picture}
\end{center}
Evaluating two morphisms at $[1_{(n]},r]_\sim \in \left[\begin{array}{c} (n] \\ n  \end{array}\right] \otimes_n R_n\subseteq \cR(n]$ after some calculations, we get
\[ \coprod_Y(\xi_Y\circ \tau_{\{Y\}}) \circ \lambda^\cR_F\circ \cR(\{-\}) ([1_{(n]},r]_\sim) = (n] \]
and
\[ (\coprod_Y\xi_Y) \circ \lambda^\cV_F\circ \tau_{\cP(n]} \circ \cR(\{-\})([1_{(n]},r]_\sim) =
(\coprod_Y\xi_Y) \circ \lambda^\cV_F \circ \cV(\{-\})([\vec{x},v]_\sim) = X \]
Thus $X=(n]$ and hance $n=m$, as required.

The naturality of $\sigma:R\ra V$ we leave for the reader. We end by showing that $\hat{\sigma}=\tau$.

By definition of $\sigma$, we have for $n\in\o$ and $r\in R_n$
\[ \tau_{(n]}([1_{(n]},r]_\sim) = [1_{(n]},\sigma_n(r)]_\sim =\hat{\sigma}([1_{(n]},r]_\sim) \]
Thus for any $ [\vec{z},r]_\sim\in  \left[\begin{array}{c} Z \\ n  \end{array}\right] \otimes_n R_n\subseteq \cR(Z)$ we have
\[ \tau_Z([\vec{z},r]_\sim) = \tau_Z(\cR(\vec{z})[1_{(n]},r]_\sim) = \]
\[ = \cV(\vec{z})(\tau_{(n]}([1_{(n]},r]_\sim)) = \cV(\vec{z})(\hat{\sigma}([1_{(n]},r]_\sim)) = \]
\[ = \hat{\sigma}_Z(\cR(\vec{z})([1_{(n]},r]_\sim)) = \hat{\sigma}_Z([\vec{z},r]_\sim) \]
i.e. $\tau=\hat{\sigma}$.  $\boxempty$
\vskip 2mm

{\em Remarks.}

1.If $(A,\alpha)$ is an $\cL$-algebra, i.e. a sup-semilattice, then the posetal collapse of the category $CP_r(A,\alpha)$ is $(A,\leq)$ i.e. the partial order on the set $A$ determined by the sup-lattice structure $(A,\alpha)$. Clearly, we have a canonical functor $CP_r(A,\alpha)\ra (A,\leq)$. The original definition of P\l onka sum on $\cR$-algebras is defined on functors $F: CP_r(A,\alpha)\ra \cEM(\cR)$ that factorize  through $(A,\leq)$.

2. The inclusion of the category of semi-analytic monads into finitary monads on $Set$ has a right adjoint. This means that any finitary monad $\cM$ on $Set$ has its 'regular part' given by
\[reg(\cM)(X)=\sum_{n\in\o} \left[\begin{array}{c} X \\ n  \end{array}\right]\otimes_n \cM(n]\]
Thus we have a canonical functor $\cEM(\cM)\lra \cEM(reg(\cM))$ induced by $reg(\cM) \ra \cM$. Hence any functor $F: C\ra \cEM(\cM)$ extends to a functor $F: C\ra \cEM(reg(\cM))$. Thus we can compute P\l onka sum of algebras for any monad except that the resulting algebra will satisfy only the regular equations of the theory.

\subsection*{Infinitary case}
The terminal monads in $\SANMND$ is the powerset monad, i.e. the monad for suplattices. We can redo the whole story above in the infinitary case getting the analogous results. We just note for the record one
\begin{theorem}\label{thm-characterization-semi-cart-infty}  Let $\tau:\cR\ra \cV$ be an arbitrary lax morphism of monads between monads in $\SANMND$. Then $\tau$ is semi-cartesian iff the induced functor $\cEM(\tau):\cEM(\cV)\ra\cEM(\cR)$ between categories of algebras preserves P\l onka sums. $\boxempty$
\end{theorem}
\vskip 2mm

An analogous result holds for monads in $\SanMnd_\kappa$, for $\kappa$ being regular cardinal.

\section{Category of linear polynomials over an algebra}\label{sec-lin-poly}

In this and the next sections we shall describe a parallel story to the one described in the previous two sections for analytic monads. We shall concentrate on the subtle differences pertinent in this case. As we saw in Section \ref{sec-corefl-cats} analytic monads are more specific than semi-analytic ones.
This has the following consequences:
\begin{enumerate}
  \item In this context we have not only a notion of a variable that occurs in a term as in regular case but we can now count occurrences; this is why we can and we will take more subtle category of polynomials, namely the category of linear polynomials.
  \item The terminal analytic monad is the monad for the commutative monoids; thus we shall equip (canonically) any category of algebra for an analytic monad with operations of P\l onka sums whose arities are categories of linear polynomials over commutative monoids.
  \item As the series expansion of analytic functors are simpler than the expansions of semi-analytic ones the definition of the morphism of monads (nat. transf. part is simpler) inducing P\l onka sums is simpler. But then the characterization results are analogous.
\end{enumerate}

Let $\cA=(\cA,\eta,\mu)$ be an analytic monad, $A:\B \ra Set$ a functor such that, for set $X$
\[ \cA(X)= \sum_{n\in\o} X^ n  \otimes_n A_n \]
$A$ is the functor part of a symmetric operad, i.e. a monoid in the monoidal category $Set^\B$ with the substitution tensor.  We define a functor
\[ CP_{l}^\cA : \cEM(\cA) \lra Cat \]
associating to an $\cA$-algebra $(X,\xi:\cA(X)\ra X)$ a category of linear polynomials $CP_l^\cA(X,\xi)$ as follows. The objects of $CP_l^\cA(X,\xi)$  are elements of $X$. A morphism in $CP_l^\cA(X,\xi)$ is an equivalence class of triples
\[ [\vec{x}, i, a ]_\sim : \vec{x}(i) \ra \alpha(\vec{x},a) \]
where, for some $n\in \o$, $\vec{x}:(n]\ra X$ is a function, $i\in(n]$, $a\in A_n$. We identify triples
 \[ \lk \vec{x}\circ\sigma, i, a \rk \sim \lk \vec{x}, \sigma(i), A(\sigma)(a) \rk \]
where $\sigma\in S_n$. The identity morphism is
\[  [x, 1, \iota ]_\sim : x \ra \xi(x,\iota)=x\]
where $\iota\in A_1$ is the unit of the symmetric operad $A$. The composition is defined by the substitution of linear-regular terms into linear-regular terms (normalization is not needed). In detail, it is defined as follows
\begin{center}
\xext=2400 \yext=280
\begin{picture}(\xext,\yext)(\xoff,\yoff)
\putmorphism(0,50)(1,0)[\vec{x}(i)`{\xi(\vec{x},a)=\vec{x}'(j)}`[\vec{x}, i, a {]}_\sim]{1200}{1}b
\putmorphism(1200,50)(1,0)[\phantom{\xi(\vec{x},a)=\vec{x}'(j)}`\xi(\vec{x}',a')`[\vec{x}', j, a' {]}_\sim]{1200}{1}b
\putmorphism(0,180)(1,0)[\phantom{\vec{x}(i)}`\phantom{\alpha(\vec{x}',a)}`[x'', i'', a''  {]}_\sim]{2400}{1}a
\end{picture}
\end{center}
where $a\in A_n$, $a'\in A_n$; the function $\vec{x}''=\vec{x}'(j\backslash\vec{x}):(n+m-1]\ra X$ is defined so that
\[ \vec{x}'(j\backslash\vec{x})(l)\;\;= \;\; \left\{ \begin{array}{ll}
                \vec{x}'(l)   & \mbox{ if } 1\leq l<j \\
                \vec{x}(l-j+1)   & \mbox{ if } j\leq l <n+j \\
                \vec{x}'(l-n+1) & \mbox{ if } n+j\leq l < n+m .
                                    \end{array}
                            \right. \]
Thus $\vec{x}'(j\backslash\vec{x})$ replaces $j$ in the domain of $\vec{x}'$ by the whole function $\vec{x}$. $a''=((\iota,\ldots,a,\ldots,\iota)\ast a')$, i.e. the composition in the symmetric operad $A$ of $a'$ with $a$ placed into the $j$'s place. This ends the definition of the category $CP_l^\cA(X,\xi)$.

A homomorphism $h:(X,\xi)\lra (X',\xi')$ induces a functor \[ CP_l^\cA(h) :CP_l^\cA(X,\xi)\lra CP_l^\cA(X',\xi')\]
so that the morphism $[\vec{x}, i, a ]_\sim : \vec{x}(i) \ra \xi(\vec{x},a)$ is sent to
\[ [h\circ\vec{x}', i, a ]_\sim : h(\vec{x}'(i)) \ra h(\xi(\vec{x},a))= \alpha'(\vec{x}',a) \]

We note for the record

\begin{fact}\label{lin-polynomials} For any analytic monad $\cA$,
the functor
\[ CP_l^\cR : \cEM(\cA) \lra Cat \]
associating to $\cA$-algebras their categories of linear polynomials is well defined. $\boxempty$
\end{fact}

\vskip 2mm
{\em Remarks.}

One can describe the category of polynomials $CP_l^\cR(X,\xi)$ as follows. Its objects are elements of the algebra $(X,\xi)$.
A morphism from $x$ to $y$ can be described in a similar way as in the semi-analytic case as follows:
\begin{enumerate}
  \item a linear-regular term in $n$ variable $a(x_1,\ldots,x_n)$ (all the variables necessarily explicitly occur in $a$ exactly once);
  \item a certain variable number $i\in(n]$ is chosen;
  \item a function $\vec{x}:(n]\ra X$ interpreting variables occurring in $a(x_1,\ldots,x_n)$;
  \item the domain on the morphism $x$ is equal to the interpretation of the $i$-th variable i.e. $\vec{x}(i)$;
  \item the codomain of the morphism $y$ is equal to the value of the term $a(x_1,\ldots,x_n)$ under the interpretation of the variables $\vec{x}$.
\end{enumerate}
The composition of such morphisms is defined as (the most reasonable) substitution.

The reason why this definition works for analytic monads is, as we explained in the introduction, that in the corresponding linear-regular equational theories there is a good notion of the number of occurrences of a variable in a term.

Note that the category of  linear polynomials of a monoid, as well as the category of regular polynomials of sup-semilattices, have natural structure of a monoidal category with multiplication (or sup) playing the role of the tensor.

\vskip 2mm

If $\tau: \cA\ra \cA'$ is a morphism of analytic monads then the functors associating categories of linear polynomials to those monads are related as follows. The morphism $\tau$ induces a 'forgetful' functor $\cEM(\tau): \cEM(\cA')\ra \cEM(\cA)$ and we have a natural transformation
\[  \gamma^\tau:  CP^\cA_l\circ \cEM(\tau) \lra CP^{\cA'}_l : \cEM(\cA')\lra Cat\]
defined for a $\cA'$-algebra $(X,\xi)$ a functor
\[ \gamma^\tau_{(X,\xi)}: CP^\cA_l(X,\xi\circ\tau) \lra CP^{\cA'}_l(X,\xi) \]
constant on object, and sending morphism $[\vec{x}, i, a ]_\sim :\vec{x}(i)\lra \xi\circ\tau(\vec{x},a)$ to the morphism
\[ [\vec{x}, i, \tau_n(a) ]_\sim : \vec{x}(i)\lra \xi(\vec{x},\tau_n(a'))  \]
where $a\in A_n$ and $\tau_n: A_n\ra A'_n$.

%\newpage
\section{P\l onka sums of algebras for an analytic monad}\label{sec-an-plonka-sums}

The theory of P\l onka sums for the categories of algebras for analytic monads is parallel to the theory of P\l onka sums for the categories of algebras for semi-analytic monads described in Section \ref{sec-reg-plonka_sum}. We need to replace categories of regular polynomials by the categories of linear polynomials and the terminal semi-analytic monad of sup-semilattices by the terminal analytic monad of commutative monoids. As the details are very similar in both cases we shall give the basic definition and state some of the facts without proofs leaving the verification as an exercise.

Let $\pi : \cA\ra \cT$ be a morphism of analytic monads, defined by a natural transformation (denoted by the same letter) $\tau:A\ra T$ in $Set^\B$. Let $(Z,\zeta)$ a $\cT$-algebra. Let denote  the category of the linear polynomials ${CP_l(Z,\zeta)}$ as $\bZ$. The monad $\hat{\cA}$ is a lift of the monad $\cA$ to the category $Set^{\bZ}$. We shall define a lax morphism of monads induced by $\pi$  and $\cR$-algebra $(Z,\zeta)$
\[ (\bigsqcup{\!_{(Z,\zeta)}},\lambda^{\cA,(Z,\zeta)}) :  \hat{\cA} \ra \cA \]
We usually drop superscripts  $(Z,\zeta)$ and write $\cA$ instead of $\pi$ in $\lambda^{\pi,(Z,\zeta)}$ when it does not lead to a confusion. We write $\bigsqcup{\!_{Z}}$ rather than $\bigsqcup{\!_{(Z,\zeta)}}$. Let $F:\bZ\ra Set$ be a functor. We shall define
\[ \lambda^\cA_F: \cA(\coprod_{z\in Z}F(z)) \lra \coprod_{z\in Z}\cA(F(z)) \]
as
\[ [\vec{x},a]_\sim \mapsto [f\circ \vec{x}',a]_\sim\in \cA(F(b)) \]
where $\vec{x}'$ is a lift of $\vec{x}$ as before, the function $f$ and element $b\in Z$ will be described below and are displayed in the diagram below. Since $[\vec{x},a]_\sim \in \cA(\coprod_{z\in Z}F(z))$, it follows that for some $n\in\o$, $\vec{x}:(n]\ra \coprod_{z\in Z} F(z)$ is a function and $a\in A_n$. Let $p:\coprod_{z\in Z} F(z)\ra Z$ be the projection from the coproduct to the index set. We put
\[ b=\zeta([ p\circ\vec{x},\pi_n(a))]) \]
and, for $i\in (n]$, we have a morphism
\[ \psi_i = [p\circ\vec{x},i,\pi_n(a)]: \vec{x}(i)\lra b\]
in the category $\bZ$. The function $f$ is

\begin{center} \xext=3200 \yext=700
\begin{picture}(\xext,\yext)(\xoff,\yoff)
\putmorphism(0,550)(1,0)[A`\coprod_{a\in A} F(a)`p]{600}{-1}a
\putmorphism(600,550)(1,0)[\phantom{\coprod_{a\in A} F(a)}`\coprod_{i\in (n]} F(\vec{a}(s(i)))`{[\kappa_{\vec{a}(s(i))}]}_{i\in(n]}]{1300}{-1}a
\putmorphism(1900,550)(1,0)[\phantom{\coprod_{i\in (n]} F(\vec{a}(s(i)))}`F(b)`f=[F(\psi_i){]}_{i\in(n]}]{1300}{1}a

          \put(1130,140){\vector(-1,1){320}}
          \put(1240,140){\vector(1,1){320}}

\put(900,200){${\vec{x}}$}
\put(1400,200){${\vec{x}'}$}
\put(1130,40){${(n]}$}
 \end{picture}
\end{center}

A simple verification shows

\begin{proposition}\label{preservation_polonka1}
$(\bigsqcup{\!_{A}},\lambda^\cA) :  \hat{\cA} \ra \cA$ is a lax morphism of monads. $\boxempty$
\end{proposition}

Lifting the above morphism of monads, we obtain  $\cT$-indexed P\l onka sum of $\cR$-algebras, i.e. for any $\cT$-algebra $(A,\alpha)$ we obtain an operation
\begin{center} \xext=1350 \yext=750
\begin{picture}(\xext,\yext)(\xoff,\yoff)
\setsqparms[1`1`1`1;1200`550]
 \putsquare(150,50)[\cEM(\cR)^{CP_l(A,\alpha)}\cong \cEM(\hat{\cR})`\cEM(\cR)`Set^{CP_l(A,\alpha)}`Set;\bigsqcup{\!_{A}}`\cU^{CP_l(A,\alpha)}`\cU`\bigsqcup{\!_{A}}]
\end{picture}
\end{center}
\vskip 2mm
The analog of Proposition \ref{T-morphisms-pres-plonka} holds verbatim for the lax morphisms of analytic monads. Thus again in this case, whenever we have $\pi : \cA\ra \cT$ a lax morphism of analytic monads, we have operations of P\l onka sums on the category of algebras for $\cA$ whose arities are the categories of linear polynomials over algebras for $\cT$.  As the category of analytic monads has the terminal object, the monad $\cC$ for commutative monoids, we can make these sums independent of a varying analytic monad $\cT$. We equip any category of algebras for an analytic monad $\cA$ with a canonical system of P\l onka sums with arities being the categories of linear polynomials over commutative monoids (i.e. algebras for the terminal analytic monad $\cC$). We have

\begin{theorem}\label{thm-characterization-weakly-cart} Let $\tau:\cA\ra \cA'$ be an arbitrary lax morphism of monads between analytic monads. Then $\tau$ is weakly cartesian iff the induced functor $\cEM(\tau)$ between categories of algebras preserves P\l onka sums. $\boxempty$
\end{theorem}

%\newpage
\section{Some examples}\label{sec-examples}
\begin{enumerate}
  \item Let $\tau: \cR\ra \cT$ be a morphism of semi-analytic monads.  $N$ is a $\cR$-algebra,  $M$ is a $\cT$-algebra. Then P\l onka sum is the usual binary product of algebras
  $$\bigsqcup_{c\in C}\bF(c)\cong N\times Alg(\tau)(M)$$
  where $\bF: C=CP_r(M) \lra Alg(S)$ is a constant functor with value $N$. For analytic monads the analogous fact holds true.
  \item Let $h:M\ra N$ be a homomorphism of $\cR$-algebras for a semi-analytic monad $\cR$. Then on the coproduct in $Set$ of universes of $M$ and $N$ there is a structure $\cR$-algebra, the P\l onka sum over the $2$-element sup-semilattice. The constants are interpreted in $M$. The operations are interpreted in $M$ if all arguments are in $M$, and in $N$ after transferring the necessary arguments from $M$ to $N$ by $h$, otherwise.
  \item The following is a more elaborate example of a P\l onka sum. It is inspired by the conversations I had with F.W. Lawvere after my talk in Coimbra during the Workshop on Category Theory in honor of George Janelidze, on the occasion of his 60th birthday.
  \begin{enumerate}
    \item Theory $Mat_R$ is the theory whose  operations from $n$ to $m$ are $m\times n$-matrices of elements of a rig $R$. Neither theory of rigs nor $Mat_R$ are regular theories, however the affine part of $Mat_R$ is a regular theory of convexity algebras denoted $CV$. It has one constant say $\perp$ and those operations $\lk r_1,\ldots, r_n\rk : n\ra 1$ that $\sum_{i=1}^n r_i = 1$, and $r_i\neq 0$ for $i=1,\ldots, n$.

    \item Now let us fix $R$ to be the rig of non-negative real numbers. Free algebra for $CV$ on $n$ generators is an $n$-dimensional simplex (with one distinguished vertex $\perp$).
    \item  Fix any set $Z$. Let $\Omega$ denote the category of regular polynomials over the sup-semilattice of finite subsets $\cP_{<\o}(Z)$. There is an obvious forgetful functor $U:\Omega\ra Set$ sending the morphism $[i, \sum_{j=1}^n r_j x_j, X_1,\ldots, X_n]: X_i\ra \bigcup_jX_j$ in $\Omega$ to the inclusion $X_i\ra \bigcup_jX_j$. The functor $\bF : \Omega \ra Alg(CV)$ is the composition of $U$ with the free $CV$-algebra functor $\cF : Set \ra Alg(CV)$.
    \item The P\l onka sum $\bigsqcup_\Omega \bF$ is the disjoint sum of simplices spanned by $\perp$ and finite subsets of $Z$. The $k$-ary operations
    $O_{\vec{r}}$ given by an $k$-tuple $\vec{r}=\lk r_1,\ldots, r_k\rk$ such that $r_i\neq 0$ and $\sum_i^kr_i=1$ in $\bigsqcup_\Omega \bF$ is defined as follows. Let $y_i\in \Delta_{X_i}$ be elements of $\bigsqcup_\Omega \bF$, i.e.  $y_i = \sum_{j\in X_i\cup\{\perp\}} s^i_j\cdot x^i_j$ for $i=1,\ldots, k$, and some $s_j^i\in R$. Then
    \[ O_{\vec{r}}(y_1,\dots, y_k) = \sum_{i=1}^k r_i \sum_{j\in X_i\cup\{\perp\}} r^i_j\cdot x^i_j \]
    and it is clearly an element of $\Delta_{\bigcup_i X_i}$.
  \end{enumerate}
  \item There are some similarities between P\l onka sums and graded rings. However there are differences, as well. A graded ring is build from a (lax) monoidal functor from a monoid $M$ considered as a discrete monoidal category to abelian groups. For example the graded ring of polynomials in $n$ variables $\mathds{Z}[x_1,\ldots,x_1]$ arises in this way from the functor $P:\o\ra Ab$ sending $n$ to the set of uniform polynomials of degree $n$, say $P_n$. The coherence transformations for this monoidal functor are multiplications $\varphi_{n,m}:P_n\times P_m\lra P_{n+m}$ and identity $\bar{\varphi}:\mathds{Z}\ra \mathds{Z}=P_0$.
\end{enumerate}

\section{P\l onka products}\label{sec-plonka-prod}

The P\l onka products are not as common as P쓾nka sums but they exists in the literature in a different setup. We describe this briefly below.

Let $(\cS,\eta,\mu)$ be a monad on $Set$, $(\hat{\cS},\hat{\eta},\hat{\mu})$ be a lift of this monad on $Set^2$, and $\mathds{1}$ be the unique monad on the terminal category $\jeden\cong Set^0$. Clearly the category $Set$ with binary products and the terminal object is a monoidal category $(Set,\times, 1,\alpha,\lambda,\varrho)$. Let $\varphi: \times \circ  \hat{\cS}\lra \cS\circ \times$ and $\bar{\varphi}: 1\ra \cS(1)$ be two natural transformations.

Recall that ${\bf MND}_{lax}$ denotes the $2$-category of monads, lax morphisms of monads, and transformations of lax morphisms of monads. This category has finite products.

Then we have

\begin{proposition}\label{polonka products}
With the notation as above
\begin{enumerate}
  \item If $(\cS,\varphi,\bar{\varphi}, \eta,\mu)$ is a monoidal monad on $(Set,\times, 1,\alpha,\lambda,\varrho)$ then
  \[ (\times, \phi) : \hat{\cS}\ra \cS\hskip 5mm {\rm and} \hskip 5mm(1,\bar{\varphi}):\mathds{1}\ra \cS\]
  are oplax morphism of monads and hence they induce P\l onka products
  \[ \bigsqcap_I:\cK(\cS)^2\ra \cK(\cS) \hskip  5mm {\rm and} \hskip 5mm \bigsqcap_\emptyset:\jeden\ra \cK(\cS)\]
  \item On the other hand, if $(\times, \phi)$ and $(1,\bar{\varphi})$ are oplax morphism of monads so that together with $(\alpha,\lambda,\varrho)$ constitute a monoidal category object on $(\cS,\eta,\mu)$ in ${\bf MND}_{lax}$ then $(\cS,\varphi,\bar{\varphi}, \mu, \eta,\mu)$ is an oplax monoidal monad on $(Set,\times, 1,\alpha,\lambda,\varrho)$ and $\cK(\cS)$ is the Kleisli object for this monad in ${\bf MND}_{lax}$. $\boxempty$
\end{enumerate}
\end{proposition}

{\em Remarks.}
From the above proposition follows that if $\cS$ is a commutative monad (c.f \cite{K}, \cite{CJ}) then, for any finite set $I$, we have an oplax morphism of monads $(\bigsqcap_I,\psi): \cS^I\ra \cS$ (where $\cS^I$ is a lift of $\cS$ to $Set^I$) and a P\l onka product $\cK(S)^I\ra \cK(S)$ ($\hat{\cS}$ is the lift of $\cS$ to $Set^I$).

Note that the (finitary) commutative monads on $Set$ form a full reflective subcategory of $\Mnd$.

%\newpage
\section{Distributive laws}\label{sec-dist-laws}
There is yet another property that the terminal monads in $\AMnd$ and $\SanMnd$ have in common: they distribute over all the other monads in `their' respective categories.

The sup-semilattice monad $\cL$ lifts canonically to the category of algebras for any regular monad due to the fact that it distributes over any semi-analytic monad  $\cR$ and the commutative monoid monad $\cC$ lifts canonically to the category of algebras for any analytic monad due to the fact that $\cC$ distributes over any analytic monad $\cA$. The former statement belongs to the folklore and the latter might be new. We shall give precise formulas describing these distributive laws.

The theory of sup-lattices has exactly one regular operation of every arity $n\in\o$.  To fix the notation we denote such an $n$-ary operation as $\mu_n$. The coefficient functor $L:\S \ra Set$ of the corresponding semi-analytic monad $\cL$ is the terminal object $Set^\S$, so that $L(n)=\{\mu_n\}$. $\cL$ is the value of the right adjoint to the inclusion functor $\SanMnd\ra \Mnd$ on the terminal monad in $\Mnd$. If $\cR$ is another semi-analytic monad with the coefficient functor $R:\S\ra Set$, then we define the distributive law of $\cL$ over $\cR$
\[ \rho_\cR : \cR\circ\cL \lra \cL \circ \cR \]
as follows. Let $RL, LR : \S\ra Set$ be the coefficient functors of the composition functors $\cR\circ\cL$ and $\cL \circ \cR$, respectively. Fix a set $X$ and an element $[\vec{x};\phi;\mu_{m_1},\ldots,\mu_{m_k};r]_\sim$ in $\left[\begin{array}{c} X \\ n  \end{array}\right] \otimes_n RL_n$, where $\vec{x}:(n]\ra X$ is an injection, $\phi:(m]\ra (n]$ is a surjection, $m=\sum_{i=1}^k m_i$ and $r\in R_k$. We put
\[ \rho_{\cR, X}([\vec{x};\phi;\mu_{m_1},\ldots,\mu_{m_k};r]_\sim) = [\vec{y};\psi;r,\ldots,r;\mu_M]_\sim \;\; \in \left[\begin{array}{c} X \\ k\cdot M \end{array}\right] \otimes_{k\cdot M} LR_{k\cdot M }\]
where $r$ occurs on the right side $M=m_1\cdot m_2\cdot \ldots\cdot m_k$ times, $\psi:(k\cdot M] \ra (n']$, $\vec{y}:(n']\ra X$ is a surjective-injective factorization of the function
\[ f: (k\cdot M] \cong (k]\times (m_1]\times \ldots \times (m_k] \lra X \]
such that for $\lk i,j_1,\ldots, j_k\rk\in (k]\times (m_1]\times \ldots \times (m_k]$
\[ f(i,j_1,\ldots, j_k) = \phi_i(j_i)\]
and $\phi_i (m_i]\ra X$ is the composition of $\phi$ with the obvious inclusion $(m_i]\ra (m]$.

Since we allow in the above definition $m_i$ to be equal $0$ so $M$ can be equal $0$ independently of the values of other $m_i$'s. This is why $\rho_\cR$ does not need to be semi-analytic natural transformation and hence $\cL\circ \cR$ does not need to be a semi-analytic monad. The situation can be improved from that point of view if we drop the constant in $\cL$, i.e. if we look at the (regular) theory of sup-semilattices without bottom element i.e. the monad $\cL'$ whose coefficient functor $L':\S\ra Set$ is as $L$ except that $L'(\emptyset)=\emptyset$. Then the above formulas do define a distributive law
\[ \rho'_\cR : \cR\circ\cL' \lra \cL' \circ \cR \]
which is semi-analytic and hence the monad $\cL'\circ \cR$ is semi-analytic. Note that $\cL'$ is the value of the right adjoint to the inclusion functor $\SanMnd\ra \Mnd$ on the only proper submonad of the terminal monad in $\Mnd$, corresponding to the theory without any operations and having one equation $x_1=x_2$.

As we mentioned above, the commutative monoid monad $\cC$ distributes over any analytic monad $\cA$. Again the theory of commutative monoids is the theory that has a unique analytic (this time) operation $\mu_n$ of arity $n$, for every $n\in\o$. The coefficient functor $C:\B\ra Set$ is the terminal object in $Set^\B$. The distributive law
\[  \alpha_\cA: \cA \circ \cC \lra \cC\circ \cA \]
is even simpler to define then $\rho_\cR$ for semi-analytic monads. Fix a set $X$. For $[\vec{x};\mu_{n_1},\ldots,\mu_{n_k};a]_\sim$ in $X^n\otimes_n AC_n$ we put
\[ \alpha_{\cA,X}([\vec{x};\mu_{n_1},\ldots,\mu_{n_k};a]_\sim) = [\vec{y};a,\ldots,a;\mu_{N}]_\sim\]
where $n=\sum_{i=1}^k n_i$, $a$ occurs on the right side $N=n_1\cdot n_2\cdot \ldots\cdot n_k$ times, and the function $\vec{y}: (k\cdot N]\ra X$ is defined for
$\lk i,j_1,\ldots, j_k\rk\in (k]\times (n_1]\times \ldots \times (n_k]$
\[ \vec{y}(i,j_1,\ldots, j_k) = \vec{x_i}(j_i)\]
and $\vec{x}_i: (n_i]\ra X$ is the composition of $\vec{x}$ with the obvious inclusion $(n_i]\ra (n]$. Again the distributive law $\alpha_\cA$ is neither analytic nor even semi-analytic. If we consider the submonad of $\cC$ without the constant i.e. the monad for the commutative semigroups $\cC'$, then we shall get a semi-analytic distributive law
\[  \alpha'_\cA: \cA \circ \cC' \lra \cC'\circ \cA \]
given by essentially the same formulas as $\alpha_\cA$.  Thus the monad $\cC'\circ \cA$ is semi-analytic but not analytic in general.

The verification of the above statements is easy but a bit tedious. We state it for the record.

\begin{theorem}\label{dist_laws-an-san}
\begin{enumerate}
  \item The monads $\cL$ and $\cL'$ for sup-semilattices and for sup-semilattices without bottom distributes over any other semi-analytic monad $\cR$ in a canonical way through the distributive laws $\rho_\cR$ and $\rho'_\cR$, respectively.
The distributive law $\rho'_\cR : \cR\circ\cL' \lra \cL' \circ \cR$ is semi-analytic and the composed monad $\cL'\circ \cR$ is again semi-analytic.
  \item The monads $\cC$ and $\cC'$ for commutative monoids and for commutative semigroups distribute over any other analytic monad $\cA$ in a canonical way through the distributive laws $\alpha_\cR$ and $\alpha'_\cA$, respectively..
The distributive law $\alpha'_\cA : \cA\circ\cC' \lra \cC' \circ \cA$ is semi-analytic and the composed monad $\cC'\circ \cA$ is semi-analytic. $\boxempty$
\end{enumerate}
\end{theorem}
The extension of these result to the infinitary cases we leave for the reader.

\end{document}